\documentclass[12pt]{amsart}
\usepackage{color}
\usepackage{amssymb}
\usepackage{graphicx}  

\newtheorem{prop}{Proposition}[section]
\newtheorem{lem}[prop]{Lemma}
\newtheorem{corollary}[prop]{Corollary}

\theoremstyle{remark}

\numberwithin{equation}{section}

\def\squarebox#1{\hbox to #1{\hfill\vbox to #1{\vfill}}}

\newcommand{\stopthm}{\hfill\hfill\vbox{\hrule\hbox{\vrule\squarebox
                 {.667em}\vrule}\hrule}\smallskip}

\def\bbbone{{\mathchoice {1\mskip-4mu {\rm{l}}} {1\mskip-4mu {\rm{l}}}
{ 1\mskip-4.5mu {\rm{l}}} { 1\mskip-5mu {\rm{l}}}}}

\newcommand{\cS}{{\mathcal S}}
\newcommand{\ra}{\rangle}
\newcommand{\la}{\langle}
\def\tn{{\mathbb T}^{2n}}
\newcommand{\F}{{\mathcal F}}

\newcommand{\Hh}{{\mathcal H}}
\newcommand{\set}[1]{\left\{\,#1\,\right\}}

\newcommand{\cL}{{\mathcal L}}

\newcommand{\tm}{{\tilde m}}
\newcommand{\tx}{{\tilde x}}
\newcommand{\hh}{{\tilde h}}
\newcommand{\txi}{{\tilde \xi}}

\renewcommand{\Re}{\mathop{\rm Re}\nolimits}
\renewcommand{\Im}{\mathop{\rm Im}\nolimits}
\newcommand{\Spec}{\operatorname{Spec}}
\newcommand{\rank}{\operatorname{rank}}
\newcommand{\tr}{\operatorname{tr}}
\newcommand{\vol}{\operatorname{vol}}
\newcommand{\supp}{\operatorname{supp}}
\newcommand{\CI}{{\mathcal C}^\infty}
\newcommand{\CIc}{{\mathcal C}^\infty_{\rm{c}}}
\newcommand{\EE}{\mathbb E}
\newcommand{\RR}{\mathbb R}
\newcommand{\TT}{\mathbb T}
\newcommand{\Sph}{\mathbb S}
\newcommand{\CC}{\mathbb C}
\newcommand{\ZZ}{\mathbb Z}
\newcommand{\NN}{\mathbb N}
\newcommand{\defeq}{\stackrel{\rm{def}}{=}}

\def \rest {\!\!\upharpoonright}
\def \restrict {\!\!\upharpoonright}
\def \mfm {{\mathfrak M}}
\def \mfp {{\mathfrak P}}
\def \mfg {{\mathfrak G}}

\ifx\pdfoutput\undefined
  \DeclareGraphicsExtensions{.pstex, .eps}
\else
  \ifx\pdfoutput\relax
    \DeclareGraphicsExtensions{.pstex, .eps}
  \else
    \ifnum\pdfoutput>0
      \DeclareGraphicsExtensions{.pdf}
    \else
      \DeclareGraphicsExtensions{.pstex, .eps}
    \fi
  \fi
\fi

\title[Probabilistic Weyl laws for quantized tori]
{Probabilistic Weyl laws for quantized tori}

   \author { T.J. Christiansen}
\address{Department of Mathematics,
University of Missouri,
Columbia, Missouri 65211, USA} 
\email{christiansent@missouri.edu}
   \author { M. Zworski}  
\address{Mathematics Department, University of
California, Evans Hall, Berkeley, CA 94720, USA}
\email{zworski@math.berkeley.edu}

\begin{document} 

\begin{abstract} For the Toeplitz quantization of complex-valued
functions
on a $2n$-dimensional torus we prove that the expected
number of eigenvalues of small random perturbations of
a quantized observable satisfies a natural Weyl law \eqref{eq:tm}.
In numerical experiments 
the same Weyl law 
also holds for ``false'' eigenvalues created by pseudospectral effects.
\end{abstract}

\maketitle

\section{Introduction and statement of the result}
\label{int}

In a series of recent papers Hager-Sj\"ostrand \cite{h-sj}, 
Sj\"ostrand \cite{Sj}, and Bordeaux Montrieux-Sj\"ostrand \cite{BoSj}
established almost sure Weyl asymptotics for small random perturbations
of non-self-adjoint elliptic operators in semiclassical and high energy
r\'egimes. The purpose of this article is to present a related simpler
result in a simpler setting of Toeplitz quantization. 
Our approach is also different: we estimate the
counting function of eigenvalues using traces rather than by studying zeros
of determinants. As in \cite{h-sj}  the singular value decomposition 
and some slightly exotic symbol classes
play a crucial r\^ole.

Thus we consider a quantization $ \CI ( \TT^{2n} ) \ni f \longmapsto
f_N \in M_{N^n} ( \CC )  $, where 
$ \TT^{2n} $ is a $2n$-dimensional torus, $ \RR^{2n}/\ZZ^{2n} $, and 
$ M_{N^n} ( \CC ) $ are $ N^n \times N^n $ complex matrices.
The general procedure will be described in \S \ref{qot} but if $ n = 1$
and $ \TT = \Sph_x \times \Sph_\xi $, then 
\begin{gather}
\label{eq:simq}
\begin{gathered}
  f = f(x) \longmapsto f_N \defeq {\rm{diag}}\, \left(f(\ell/N)
\right)  \,, \ \ {\ell=0,\cdots , N-1}\,,  \\
g = g( \xi ) \longmapsto g_N \defeq {\mathcal F}_N^* \,  {\rm{diag}}\,
\left( g(k/N) \right) {\mathcal F}_N \,,  \ \ 
{k=0,\cdots , N-1} \,, 
\end{gathered}
\end{gather}
where $ {\mathcal F}_N = ( \exp ( 2 \pi i k \ell/ N) /
\sqrt {N } )_{0\leq k, \ell \leq 0 , N-1 } $, is the discrete Fourier transform.

Let $ \omega \mapsto Q_N ( \omega ) $ be the gaussian ensemble of 
complex random $ N^n \times N^n $ matrices -- see \S \ref{rmt}. 
With this notation in place we can state our result:

\noindent
{\bf Theorem.}{\em \,\,Suppose that $ f \in \CI ( \TT^{2n} )  $, and that 
$\Omega$ is a simply connected 
open set with a smooth boundary, $ \partial \Omega $, such that
for all $ z  $ a neighbourhood of $ \partial \Omega $, 
\begin{equation}   {\rm{vol}}_{{\mathbb T}^{2n}} \left( \{ w \; : \;
| f ( w ) - z | \leq t \} \right) = {\mathcal O} ( t^{\kappa } ) \,, \ \
0 \leq t \ll 1 \,, \label{eq:cond}
\end{equation}
with  
$  1 < \kappa \leq 2 $. 
Then
 for any $ p \geq p_0 > n+1/2$ 
\begin{equation}
\label{eq:tm}
\EE_\omega \left( | \Spec ( f_N + N^{-p} Q_N (\omega ) ) \cap 
\Omega | \right) = N^n \vol_{\TT^{2n}} ( f^{-1} ( \Omega ) ) + 
{\mathcal O}  ( N^{n-\beta} ) \,, \end{equation}
for any $  \beta < ( \kappa - 1 ) / ( \kappa + 1 ) $.
}

\begin{figure}
\begin{center}
\includegraphics[width=6.3in]{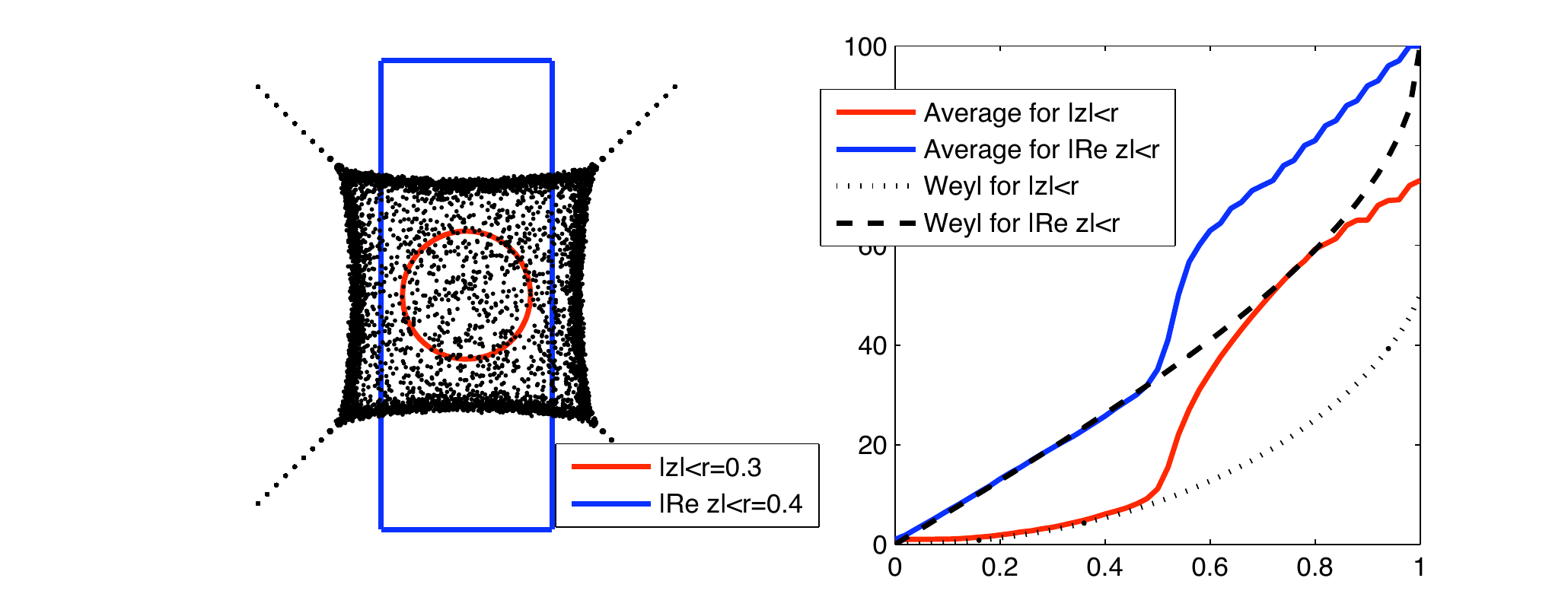}
\end{center}
\caption{On the left we reproduce \cite[Fig.1.2]{ChTr} with two 
types of regions used for counting added. It represents 
$ \Spec ( f_N + E ) $ where $ N = 100 $, 
$ f ( x , \xi ) = \cos ( 2 \pi x ) + i \cos ( 2 \pi \xi ) $ 
(called the ``the Scottish flag operator'' in \cite{ChTr}), for 
a hundred complex random matrices, $ E $, of norm $ 10^{-4} $. On the right
we show the counting functions for the two regions, and the corresponding
Weyl laws, as functions of $ r $. The breakdown of the Weyl law 
approximation occurs when the norm of the resolvent $ (f_N -z )^{-1} $,
$ | z | = r $, or $ | \Re z | = r $, is smaller than $ \| E \|^{-1} = 10^4 $.
For $ \Omega = \{ | z| < r \} $, $ r < 1 $, $ \kappa = 2 $ and for $ 
\Omega = \{ | \Re z | < r \} $,
$ \kappa = 3/2 $ at four points of $ \partial \Omega $ (intersection
with the boundary of $ f ( \TT ) $). For $ r = 1 $, the corners satisfy 
\eqref{eq:cond} with $ \kappa = 1$. }
\label{f:0}
\end{figure}

\begin{figure}
\begin{center}
\includegraphics[width=6in]{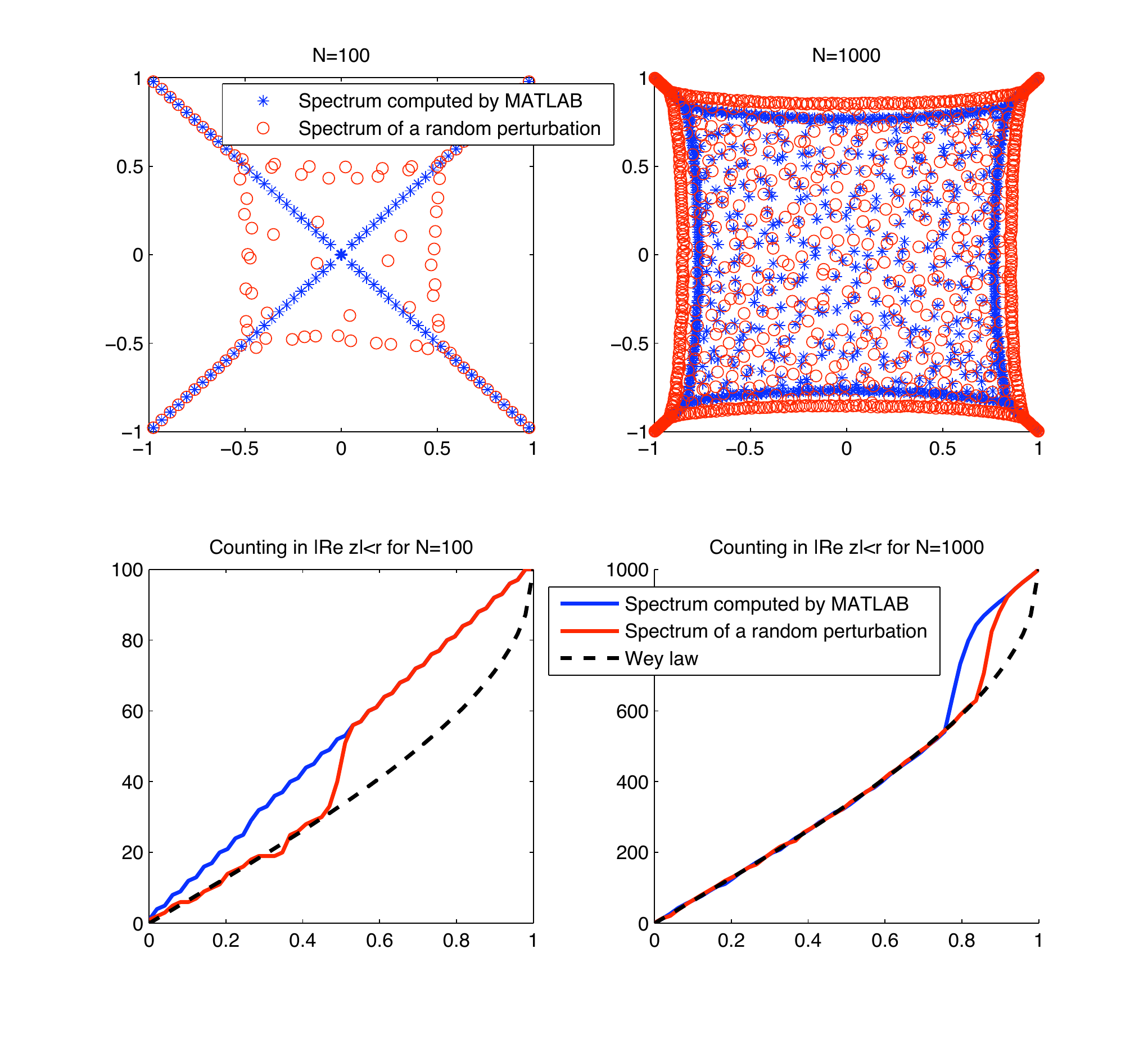}
\end{center}
\caption{The MATLAB computed spectra of $ f_N $ for 
$ f ( x , \xi ) = \cos ( 2 \pi x ) + i \cos ( 2 \pi \xi ) $. For $ N = 100$
the computations return the correct eigenvalues following the 
Scottish flag pattern. For $ N = 1000 $ the actual spectrum of $ f_N $
still follows the same pattern but the computations return ``false eigenvalues''
which satisfy the same Weyl law as random perturbations.
The plots of the counting function for a {\em single} random matrix 
are very close the Weyl asymptotic even in the case of $ N = 100$
providing support for the conjecture in \S \ref{int}.}
\label{f:2}
\end{figure}

\noindent
{\bf Remark.} The theorem applies to more general operators
of the form $ A(N) = f_N + g_N/N $, where $ g $ may depend
on $ N $ but all its derivatives are bounded as $ N \rightarrow 
\infty $.

The main point of the probabilistic Weyl law 
\eqref{eq:tm} is that for most complex-valued functions $ f $
the spectrum of $ f_N $ will {\em not} satisfy the Weyl law.
Yet, after adding a tiny random perturbation, the spectrum 
will satisfy it in a probabilistic sense. As illustrated in 
Fig.\ref{f:2} a tiny perturbation can change the spectrum
dramatically, with the density of the resulting eigenvalues 
 asymptotically determined by the original function $ f$.

Condition \eqref{eq:cond} with $ 0 < \kappa \leq 2 $ appears in the
work of Hager-Sj\"ostrand \cite{h-sj}. Its main r\^ole here is to 
control the number of small eigenvalues of $ ( f_N - z)^* ( f_N - z ) $,
see Proposition \ref{p:rank}, and that forces us to restrict to the case
$ \kappa > 1 $.
It is a form of a {\L}ojasiewicz inequality and for 
real analytic $ f $ it always holds for some $ \kappa > 0 $, as
can be deduced from a local resolution of
singularities -- see \cite[Sect.4]{BM}. Similarly, for $ f $ {\em 
real analytic} and such that $ f ( \tn ) \subset \CC $ has a non-empty
interior,
\[  df \rest_{ f^{-1} ( z ) } \neq 0 \ \Longrightarrow \ 
\text{ \eqref{eq:cond} holds with $ \kappa > 1 $.} \]
For $ f \in \CI ( \TT^{2n} ) $ we have 
\[ df \wedge d \bar f \rest_{ f^{-1} ( z ) } \neq 0 \ \Longrightarrow \ 
\text{ \eqref{eq:cond} holds with $ \kappa = 2 $,} \]
and by the Morse-Sard theorem the condition on the left 
is valid on a complement of a set
of measure $ 0 $ in $ \CC $. Also, 
\[ \forall \,  w \in f^{-1} ( z)  \ \ \{ f, \bar f \} ( w ) 
\neq 0 \ \text{or} \ \{ f , \{ f , \bar f \} \} ( w ) \neq 0 
\ \Longrightarrow \ 
\text{ \eqref{eq:cond} holds with $ \kappa = \frac 3 2  $,} \]
see \cite[Example 12.1]{h-sj}. Here $ \{ \bullet , \bullet \} $
is the Poisson bracket on $ \TT^{2n} $ (see \eqref{eq:po2co} below).

The significance of the Poisson bracket in this context comes
from the following fact:  
\begin{equation} 
\label{eq:Hoer} \{\Re f  , \Im f \} ( w ) <  0 \,, \ \ z = f ( w ) 
\ \Longrightarrow \ \| ( f_N - z )^{-1} \| > N^p/ C_p  \,, \ 
\forall \, p > 0 \,, \end{equation}
and moreover an approximate eigenvector, $ u_N$, causing the growth of the
resolvent can be microlocalized at $ w $ (meaning that 
for any $ g $ vanishing near $ w $, $ \| g_N u_N \|_{\ell^2} = 
{\mathcal O}( N^{-\infty}  ) $, $ \| u_N \|_{\ell^2} = 1 $, see \S \ref{qot}).
This is a reinterpretation of a now classical result of 
H\"ormander proved in the context of solvability of partial
differential equations -- see \cite{DeSjZw}, \cite{Zw}, 
and references given there. For quantization of $ {\mathbb T} $
\eqref{eq:Hoer} was proved in \cite{ChTr}, and for general 
Berezin-Toeplitz quantization of compact symplectic K\"ahler 
manifolds, in \cite{BoU}.

The relation \eqref{eq:Hoer} shows that
$ \{ \bar f , f \} \neq 0 $ implies the instability of the
spectrum under small perturbations. In that case 
the theorem above is most interesting, as shown in Figures \ref{f:0}
and \ref{f:2}. However, as stressed in \cite{BoSj},\cite{h-sj}, and \cite{Sj},
the results on Weyl laws for small random perturbations have
in themselves nothing to do with spectral instability. For
normal operators they do not produce new
results compared to the standard semiclassical Weyl laws:
the distribution of eigenvalues is not affected by small
perturbations and satisfies a Weyl law to start with.

The numerical experiments suggest that much stronger results
then our theorem are true. In particular we can formulate the
following

\noindent
{\bf Conjecture.} {\em Suppose that \eqref{eq:cond} holds for all
$ z \in \CC $ with a fixed $ 0 < \kappa \leq 2 $. 
Define { random} probability measures:
\[  \mu_N ( \omega ) = \frac 1 {N^n} \sum_{ z \in 
{\rm{Spec}}\,  ( f_N + N^{-p} R_N( \omega) ) }  \delta_z \,. \]
Then, {almost surely in $ \omega $},
\[  \mu_N ( \omega ) \longrightarrow f_* ( \sigma^n/ n!) \,, \ \ \ 
N \longrightarrow \infty \,, \]
where $ \sigma = \sum_{k=1}^n d\xi_k \wedge d x_k $, $ ( x, \xi) \in 
\TT^{2n} $, is the symplectic form in $ \TT^{2n}$.}

\bigskip


The result should also hold for more general ensembles than 
complex Gaussian random matrices. Sj\"ostrand's recent paper \cite{Sj}
suggests that random {\em diagonal} matrices would be enough
to produce the Weyl law-creating perturbations. 

Bordeaux Montrieux \cite{Bo} pointed out to us that 
by taking singular $ f$'s, or $ f$'s for which derivatives
grow fast in $ N $ (corresponding to $ \rho = 1 $ in the 
$ S_\rho $ classes described in \S \ref{pc}), usual Toeplitz
matrices fit in this scheme and that numerical experiments
indicate the validity of Weyl laws in this case.

Hager \cite{Ha} indicated how the methods of \cite{h-sj} should
apply to the case of Berezin-Toeplitz quantization but that
approach did not suggest any simplifications in the method.
In this paper we follow the most na\"{\i}ve approach
which starts with the following {\em false} proof of the theorem:
\[ \begin{split}
| {\rm {Spec} } \, (f_N ) \cap \Omega | & =
 \frac{1}{ 2 \pi i }   \int_{\partial \Omega }
\tr ( f_N - z )^{-1}  dz \\
&  \! {\text{``=''}} \, \, N^n  \frac{1}{ 2 \pi i } \int_{\partial \Omega} 
\int_{ \TT^{2n} } ( f ( w) - z )^{-1}
d {\mathcal L} ( w) d z + o ( N^n ) \\
&  =  N^n \int_{ \TT^{2n} }   \left(  \frac{1}{ 2 \pi i }   \int_{\partial \Omega}
( f ( w) - z )^{-1} dz  \right)
d {\mathcal L} ( w) d z + o ( N^n ) 
  \\ & = N^n {\rm{vol}}_{ \TT^{2n} } ( f^{-1} ( \Omega ) ) + o ( N^n ) \,.
\end{split} \]
Here we attempted to apply Lemma \ref{l:trace} below as if
$ ( f_N - z )^{-1} = g_N $ for some nice function $ g$. As \eqref{eq:Hoer}
shows that is impossible in general. The random perturbation,
and taking of expected values, make this argument rigorous. 
In \S \ref{rdp} we show how the first integral split to integrals 
over small (that is of size $ \sim N^{-1/2 + \epsilon } $)  subintervals
of $ \partial \Omega $ can be replaced by integrals of 
invertible operators. That is done using the singular value
decomposition (see \cite[\S 3.6]{SjZw07} for a simple related example)
and facts about random matrices proved in \S \ref{rmt}. Based 
on the material reviewed in \S \ref{qot} we further reduce the
analysis to that of traces of an inverse of an operator which
is a quantization of a slightly exotic function on the torus.
Here ``slightly exotic'' refers to the behaviour of derivatives
as $ N \rightarrow \infty $. An application of a semiclassical
calculus gives the desired trace and concludes the proof.

Except for some facts about the standard semiclassical calculus
of pseudodifferential operator recalled in \S \ref{pc}, the 
paper is meant to be self-contained. One of the advantages
of Toeplitz quantization is the ease with which traces and 
determinants can be taken, without worries associated with 
infinite dimensional spaces.

\medskip
\noindent
{\sc Acknowledgments.}
We would like to thank Edward Bierstone and Pierre Milman for 
helpful discussions of {\L}ojasiewicz inequalities,
Mark Rudelson for suggestions concerning random matrices, 
and St\'ephane Nonnenmacher for comments on early versions
of the paper.
The authors gratefully acknowledge the partial support
by an MU research leave,
and NSF grants DMS 0500267, DMS 0654436.  
The first author thanks the Mathematics 
of U.C. Berkeley for its hospitality in spring 2009.

\section{Quantization of tori}
\label{qot}

The Toeplitz quantization of tori, or of more  
general classes of compact symplectic manifolds,
has a long tradition and we refer to \cite{BouzDB} for references
in the case of tori, and to \cite{Bor} for the case of
compact symplectic K\"ahler manifolds. 
We take a lowbrow approach and our presentation which follows \cite{NZ} is self-contained but assumes
the knowledge of standard semiclassical calculus in $ \RR^n $.
It is reviewed in \S \ref{pc} with 
detailed references to \cite{DiSj} and \cite{EZ} provided. To see
how this fits in the more general scheme see for instance \cite[\S 4.2]{BoU}. 

\subsection{Review of pseudodifferential calculus in $ \RR^n $}
\label{pc}
We first recall from \cite[Chapter 7]{DiSj} (see also \cite[Chapter 3]{EZ}) 
the quantization of functions 
$ a \in  S_\rho ( T^*\RR^n )$, 
\begin{gather*}
S_\rho  ( T^* \RR^n ) \defeq
\{ a \in \CI ( T^* \RR^n ) \; : \;  \forall \alpha,\beta\in\NN^n,\;
| \partial^\alpha_x 
\partial_\xi^\beta a ( x , \xi )  | \leq C_{\alpha \beta} h^{-(|\alpha|
+ |\beta|) \rho} \} \,, \\  0 \leq \rho < \frac 12 \,. 
\end{gather*}
To any $ a \in \cS ( T^* \RR^n ) $ we associate its $h$-Weyl 
quantization,
that is the operator $a^w ( x, hD )$ acting as follows on $\psi\in\cS(\RR^n)$:
\begin{equation}
\label{eq:aw} [a^w ( x, hD )\, \psi](x)\defeq \frac{1}{( 2 \pi h )^n } 
\int \! \int a \Big( \frac{x+y}{2} , \xi \Big)\, e^{ \frac{i}{h}
\la x - y, \xi \ra}\, \psi (y)\, dy\, d\xi \,.
\end{equation}
This operator is easily seen to have the following mapping properties
$$
a^w ( x , h D ) \; : \; \cS ( \RR^n ) \; \longrightarrow 
\; \cS ( \RR^n ) \,, \ \ 
  a^w ( x , h D ) \; : \; \cS' ( \RR^n ) \; \longrightarrow 
\; \cS' ( \RR^n ) \,,
$$
see for instance \cite[\S 3.1]{EZ}  for basic properties of the Schwartz space 
$ \cS $ and \cite[\S 4.3.2]{EZ} for the mapping properties. It can then 
be shown \cite[Lemma 7.8]{DiSj} that $ a \mapsto a^w ( x , h D ) $ 
can be extended to any $ a\in S_\rho $, and that the resulting 
operator has the same mapping properties. Furthermore, $a^w ( x , h D )$ is 
a bounded operator on $L^2(\RR^n)$. 
The condition $ \rho < 1/2 $
is crucial for the asymptotic expansion in the
 the composition formula for 
pseudodifferential operators. If $ a , b \in S_\rho $ then 
\begin{gather}
\label{eq:comp}
\begin{gathered}
a^w ( x , h D ) \circ b^w ( x , h D ) = c^w ( x , h D ) \,,  \ \ c  
= a \#_h b \in S_\rho \,,
\\
c( x, \xi ) \sim \sum_{k=0}^\infty \frac 1 {k!} \left( \frac {i h }
2 \sigma ( D_z , D_w ) \right)^k a ( z ) b ( w ) \rest_{ z = w = ( x, \xi ) } \,,\\
\text{where }\sigma(z,w)=\sigma(z_1,z_2,w_1,w_2)= \langle z_2, w_1\rangle -
\langle z_1,w_2\rangle\, .
\end{gathered}
\end{gather}
We note that 
\[ \frac 1 {k!} \left( \frac {i h }
2 \sigma ( D_z , D_w ) \right)^k a ( z ) b ( w ) = {\mathcal O} ( h^{k(1-2\rho ) } )
\,, \]
so that the expansion in \eqref{eq:comp} makes sense asymptotically. 

It is important to recall the standard way in which the quantization 
of $ S_\rho (T^* \RR^n ) $ reduces to the quantization of 
\[ S ( T^* \RR^n ) 
\defeq S_0 ( T^*\RR^n ) \,, \]
with a new semiclassical parameter, 
$  \tilde h = h^{1 - 2 \rho } \,. $
Define
$ ( \tilde x , \tilde \xi ) = ( \tilde h / h )^{\frac12} ( x , \xi ) $, 
and a unitary operator on $ L^2 ( \RR^n ) $:
\[ U_{ h , \tilde h } u ( \tilde x ) = ( \tilde h / h )^{\frac{n}{4}}
u ( ( h / \tilde h )^{\frac12}  \tilde x ) \,.\]
Then 
\[ a ( x , h D_ x ) = U_{ h , \tilde h }^{-1} \tilde a 
( \tilde x , \tilde h D_{\tilde x } ) U_{ h , \tilde h} \,, \ \ 
\tilde a ( \tilde x , \tilde \xi ) \defeq a ( 
( h / \tilde h)^{\frac12} ( \tilde x , \tilde \xi ) ) \,. \]
We have 
\[ a \in S_\rho ( T^* \RR^n ) \ \Longleftrightarrow \ \tilde a \in S ( T^* 
\RR^n ) \,. \]

One simple application of this rescaling is a version 
of the semiclassical Beals Lemma \cite[Chapter 8]{DiSj} (see also 
\cite[\S 8.6]{EZ}):
\begin{gather}
\label{eq:Beals}
\begin{gathered}
A = a^w ( x, h D) \,, \ \ a \in S_\rho ( T^*\RR^n ) \
\Longleftrightarrow 
\
 {\rm{ad}}\,_{ \ell^w_1} \circ \cdots {\rm{ad}}\,_{\ell^w_N  } A = 
{\mathcal O}_{L^2 \rightarrow L^2} (h^{N(1-\rho) } ) \,, 
\\
\text{for any sequence $ \{ \ell_j \}_{j=1}^N $ of linear functions on 
$ T^* \RR^n $}
\end{gathered}
\end{gather}

The composition formula \eqref{eq:comp} holds also for operators
in more general symbol classes. For reasons  which should 
become clear below, we will discuss it only for the $\tilde h$-quantization
with $ \rho = 0 $. First we need to recall the definition of 
an order function: $ \tm = \tm ( \tx, \txi ) $ is an order function 
if there exist $ C $ and $ M $ such that for all $ ( \tx , \txi ) $ and
$ ( \tx', \txi' ) $, we have
\[ m ( \tx , \txi) \leq C m( \tx' , \txi') ( 1 + d_{\RR^{2n}} ( 
( \tx, \txi), ( \tx' , \txi')))^M \,. \]
We then say that $ \tilde a \in S ( \tm ) $ if for all $ \alpha $, 
$ |\partial^{\alpha }_{\tx , \txi}  \tilde a 
( \tx , \txi ) | \leq C_{\alpha} \tm ( \tx , \txi ) $.
If $  \tm_1 $ and $ \tm_2 $ are two order functions and $ \tilde a 
\in S ( \tm_1 ) $, $ b \in S ( \tm_2 ) $ then 
$  \tilde a ( \tx , \hh D) \circ  \tilde b ( \tx , \hh D) = 
 \tilde c  ( \tx , \hh D) $, $\tilde c \in S ( \tm_1 \tm_2 ) $, 
and the asymptotic expansion \eqref{eq:comp} is valid in $ S ( \tm_1 \tm_2 ) $.

This has a standard application which will be crucial in \S \ref{pot}:
\begin{gather}
\label{eq:elli}
\begin{gathered}
\tilde a \in S ( \tm ) \,, \  \forall \, ( \tx, \txi ) \,, \ \ 
| a ( \tx, \txi) | \geq \tm ( \tx , \txi) \ \Longrightarrow \\
\exists \, \hh_0 \, \forall \, 0 < \hh < \hh_0 \,, \ \ \tilde a^w 
(\tx , \hh D)^{-1} 
= \tilde b^w ( \tx , \hh D ) \,, \ \ b \in S ( 1/\tm ) \,, \ \ 
\end{gathered}
\end{gather}
see for instance \cite[\S 4.5,\S 8.6]{EZ}.

The reason that we presented the order functions on the $\hh$-side
is motivated by the fact that we need the rescaling of these
order functions on the $ \hh$-side: we say that $ m = m ( x , \xi ) $
is an $ h^{\rho}$-order function if 
there exist $ C $ and $ M $ such that for all $ ( x , \xi ) $ and
$ ( x', \xi' ) $, we have
\begin{equation}
\label{eq:epsor0}
 m ( x , \xi) \leq C m( x' , \xi') ( 1 + d_{\RR^{2n}} ( 
h^{-\rho} ( x, \xi),  h^{-\rho} ( x' , \xi')))^M\,, 
\end{equation}
which means that $ \tilde m ( \tx , \txi ) \defeq m ( h^\rho \tx , h^\rho \txi ) 
$ is a standard order function defined above.
The symbol class is defined analogously, $ a \in S_\rho ( m ) $ if
$ \partial^\alpha a = {\mathcal O} ( h^{-|\alpha| \rho} m ) $. By the 
rescaling argument the ellipticity statement \eqref{eq:elli} is
still applicable if $ \rho < 1/2 $. 

The following $ h^\rho $-order function coming from 
\cite[\S 4]{h-sj} will be essential to our arguments here, and 
in \S \ref{pot} (Lemma \ref{l:ord}):
\begin{lem}
\label{l:ord0}
For $ a \in S ( T^* \RR^n ) $ 
\[  m ( x , \xi ) \defeq |a ( x , \xi ) |^2 + h^{2 \rho} \,, 
\ \ \ 0 \leq \rho < \frac 12 \,, 
\]
is an $ h^{2\rho} $-order function in the sense of definition \eqref{eq:epsor0}.
In addition, 
for $  \psi \in \CIc ( \RR ; [ 0 , 1 ] ) $ equal to $ 1 $ on $ [ - 1, 1] $,
\begin{equation}
\label{eq:mell0}   \left( | a ( x, \xi  )|^2 +  h^{2\rho}  \psi \left( 
\frac {| a ( x, \xi)|^2 }{ h^{2\rho} } \right) \right)^{\pm 1} 
 \in S_\rho  ( m^{\pm 1}  ) \,. \end{equation}
\end{lem}
\begin{proof} 
This follows from the arguments in \cite[\S 4]{h-sj} but for 
the reader's convenience we present an adapted version. 
We will use the notation $ ( \tilde x, \tilde \xi) $ 
introduced above, with $ \hh = h^{ 1 - 2 \rho }$. 
Let us put $ F ( \tx, \txi)  \defeq |a  (  x , \xi ) |^2 $, 
so that $ m  ( x, \xi ) = h^{2\rho}  \tm ( \tx , \txi)$ , where 
\[ \widetilde m ( \tx , \txi) 
\defeq h^{-2\rho} F ( \tx , \txi)  +  1 \geq 1 \,. \]
To prove \eqref{eq:epsor0} we need 
\begin{equation}
\label{eq:tmor} \widetilde m ( w ) \leq C 
{ \widetilde m ( w' )
 } ( 1 +  d_{\RR^{2n} } (  w ,  w')) ^M  \,,
\end{equation}
For $ |\beta | = 1 $,  $ \partial^\beta F = {\mathcal O} ( h^{\rho} \sqrt F ) 
$ and hence 
\[  \partial^\beta \widetilde m = \frac 1 {h^{2 \rho}} \partial^\beta F 
 = 
 { \mathcal O} ( h^{-\rho} \sqrt F ) = 
 {\mathcal O} ( \sqrt {\widetilde m }) \,. \]
For $ |\beta | = 2 $, $ \partial^\beta F = {\mathcal O} ( h^{2 \rho} ) $, 
and hence $ \partial^\beta \tm = {\mathcal O} ( 1 ) $.
By Taylor's formula,
\[ \begin{split} \widetilde m (   w') 
& \leq \widetilde m (   w ) 
+ C \sqrt{ \widetilde m (   w ) } d_{ \RR^{2n} }  
(   w,   w') 
 + C  d_{ \RR^{n} }  (   w,   w')^2 \\ 
& \leq C ( 1 + \widetilde m (   w) ) ( 1 +  d_{ \RR^{2n} }  
(   w,   w'))^2 
\,. \end{split} \]
As $ \widetilde m \geq 1 $ this proves \eqref{eq:tmor} with $ M = 2$, and 
consequently the
first part of the lemma.

For the second part we observe that $ \psi ( | a |^2 / h^{2 \rho } ) \in 
S_\rho( 1 )  $ and hence $ h^{2\rho} \psi ( | a |^2 / h^{2 \rho } )  
\in S_\rho ( m ) $.  This means that we already have the $ + $ case
of \eqref{eq:mell0}. But, 
\[    | a ( x, \xi  )|^2 + h^{2\rho}  \psi \left( 
 {| a ( x, \xi)|^2 }/ { h^{2\rho} } \right)  \geq m ( x, \xi ) /2 \,,\]
and the $ - $ case follows.
\end{proof}

We remark that by introducing $ \tilde h $ as a small, eventually fixed,
parameter, we can include the case of $ \rho = 1/2 $ -- see for 
instance \cite[\S 3.3]{SjZw04}. That type of calculus is 
used in \cite{h-sj}.

The last item in this review is a slightly non-standard 
functional calculus lemma:
\begin{lem}
\label{l:func0}
Suppose that $ a \in S_0 ( T^* \RR^n )  $, $0 \leq \rho < 1 /2 $, and 
that $ \psi \in \CIc ( \RR ) $. 
Then 
\begin{gather}
\label{eq:func0}
\begin{gathered}
\psi \left( a^w ( x , h D ) a^w ( x , h D )^*   / h ^{2 \rho } \right) = q^w ( x , h D )  \,, \ \
q \in S_\rho ( T^* \RR^n )\,, \\
q = q_0 + h^{1- 2 \rho} q_1 + {\mathcal O}_S ( h^\infty) \,, 
\ \  q_j \in  S_\rho\,, \ \
q_0 ( x, \xi )  = \psi ( |a ( x , \xi) ) |^2/ h^{2\rho} ) \,, \\
q_1 ( x, \xi ) = \widetilde \psi ( |a ( x, \xi ) |^2/ h^{2 \rho} ) 
\widetilde q_1 ( x, \xi) \,, \ \ \widetilde q_1 \in S_\rho \,, \ \ 
\widetilde \psi \in \CIc ( \RR ) \,, \ \ \widetilde \psi\rest_{{\rm{supp}}\, \psi} 
\equiv 1 \,. 
\end{gathered}
\end{gather}
\end{lem}
\begin{proof}
This is a simpler version of \cite[Proposition 4.1]{h-sj} which follows
the approach to functional calculus of pseudodifferential operators
based on the Helffer-Sj\"ostrand formula for a function of a selfadjoint operator $ A $:
\begin{equation}
\label{eq:HeSj}
 \psi ( A ) = \frac{ 1} \pi \int_{\CC} ( z - A)^{-1} \partial_{\bar z} \tilde 
\psi ( z ) d{\mathcal L} ( z ) \,, \ \ \psi \in \CIc ( \RR ) \,, 
\end{equation}
where $ \tilde \psi \in \CIc ( \CC) $ is an 
almost analytic extension of $ \psi $, $ \tilde \psi
\rest_{\RR } = \psi $ and $ \partial_{\bar z } \tilde \psi = {\mathcal O } 
( |\Im z |^\infty ) $ -- see \cite[Chapter 8]{DiSj} and
references given there. The 
reduction to the case given in \cite[Theorem 8.7]{DiSj} proceeds
as follows: the operator $ a^w ( x , h D ) a^w ( x , h D )^*   / h ^{2 \rho } 
= b^w ( x , h D)$, where $ b \in S_\rho ( m_1 ) $, where 
$ m_1 $ is an $h^\rho$-order function given by $ h^{-2\rho} m $, where
$ m $ is given in Lemma \ref{l:ord0}. By the rescaling argument above,
which gives a reduction to the case of the calculus with $ \tilde h = 
h^{1 - 2 \rho } $, we can apply \cite[Theorem 8.7]{DiSj} which 
gives that $ \psi (  a^w ( x , h D ) a^w ( x , h D )^*   / h ^{2 \rho }  ) =
g^w ( x, h D ) $, where 
$ g \in S_\rho ( m_1^{-1} ) \subset S_\rho ( 1 ) $. The symbolic expansion
presented in \cite[Chapter 7]{DiSj} complete the proof.
\end{proof}

\subsection{Quantum space associated to $ \tn$}
To define this finite dimensional space we
fix our notation for the Fourier transform on ${\mathcal S}' ( \RR^n )$: 
$$
\F_h u (\xi) \defeq \frac{1}{( 2 \pi h )^{n/2} }
\int u (x)\,e^{-\frac{i}{h} \langle x , \xi \rangle }\, dx \,,  
\qquad \F_h^* = \F_h^{-1} \,, 
$$
and as usual in quantum mechanics, 
$\F_h u(\xi)$ is the ``momentum representation'' of the state $u$.
To find the space of states we consider 
distributions $u\in\cS'(\RR^n)$ which are periodic in both position
and momentum:
\begin{equation}
\label{eq:quasi}
u ( x + \ell ) = 
u ( x ) \,, \qquad 
\F_h u (\xi + \ell ) = \F_h u (\xi)\,,
\end{equation}
see \cite[\S 4.1]{NZ} and references given there for more 
general spaces with different 
Bloch angles.
Let us denote by $ \Hh_h^n $ the space of distributions satisfying 
\eqref{eq:quasi}. The following lemma is easy to prove.
\begin{lem}
\label{l:elem}
$ \Hh_h^n \neq \{ 0 \} $ if and only if $ h = ( 2 \pi N )^{-1} $ for some positive
integer $N$, in which 
case $ \dim \Hh_h^n = N^n $ and 
\begin{equation}
\label{eq:bas}  \Hh_h^n = {\rm{span}} \set{ \frac{1}{\sqrt{N^n}}
\sum_{ \ell \in \ZZ^n } 
\delta( x - \ell - j/N) \; : \; j \in (\ZZ / N \ZZ)^n } \,.
\end{equation}
\end{lem}
For $h = ( 2 \pi N )^{-1} $, the Fourier transform $\F_h$ maps
$\Hh_h^n$ to itself. In the above basis, it is represented by the discrete Fourier 
transform
\begin{equation}\label{e:DFT}
(\F_N)_{j,j'}=\frac{e^{-2 i \pi \la j,j'\ra/N}}{N^{n/2}}\,,
\quad j,j'\in (\ZZ/ N \ZZ)^n\,.
\end{equation}
The Hilbert space structure on $ \Hh_h $ will be determined (up to 
a constant) once we define the quantization procedure. That will
be done by demanding that real functions are quantized into 
self-adjoint operators. 

\subsection{Quantization of $ \CI ( \tn ) $.}
\label{qoc}
The definition \eqref{eq:aw} immediately shows that for $ f \in \CI $
satisfying
\[ \forall \, \ell, m, \in \ZZ^n\,, \ \ \ 
 f ( x + \ell , \xi + m) = f ( x , \xi ) \ \Longrightarrow \ 
 f^w ( x , h D )\; : \; \Hh_h^n \longrightarrow \Hh_h^n \,,  \]
where we consider  $\Hh_h^n \subset \cS' ( \RR^n ) $.
Identifying a function $ f \in \CI (\tn) $ with a 
periodic function on $ \RR^{2n} $, we define
\begin{gather*}
   f_N = f^w ( x, h D)\rest_{ \Hh_h^n } \,, \ \ h = \frac 1 { 2 \pi N } \,, \ 
\ \
\CI ( \tn ) \ni f \; \longmapsto \; f_N \in 
\cL (\Hh_h^n, \Hh_h^n )\,,
\end{gather*}
and we remark that $ 1_N = Id_{\Hh_h^n}  $.

The composition formula from \S \ref{pc} applies 
since $ a , b \in \CI ( \tn ) $ can be identified with 
periodic functions on $ \RR^{2n} \simeq T^* \RR^n $ and 
\begin{equation}
\label{eq:comptor} a_N \circ b_N = c_N \,, \ \ c = a \#_h b \,, \ \ h = \frac 1 {2 \pi N}
\,, \end{equation}
where $ a \#_hb $ is as in \eqref{eq:comp}. This means that 
we simply use the standard pseudodifferential calculus but act
on a very special finite dimensional space.

 The Hilbert space structure
on $ \Hh_h^n $ is determined by the following simple result 
\cite[Lemma 4.3]{NZ} which we recall below. 
\begin{lem}
\label{l:hilb}
There exists a unique (up to a multiplicative constant) Hilbert
structure on $ \Hh_h^n $ for which all $ f_N \; : \; \Hh_h^n \; 
\rightarrow \; \Hh_h^n $ with $ f\in\CI(\tn;\RR) $ are 
self-adjoint. 
One can choose the constant so that the basis 
in \eqref{eq:bas} is orthonormal. 
This implies that the Fourier transform on $\Hh_h^n$ (represented by the
unitary matrix \eqref{e:DFT}) is unitary. 
\end{lem}
\begin{proof}
Let $ \langle \bullet , \bullet \rangle_0 $ be the inner product 
for which the basis in \eqref{eq:bas} is orthonormal, and 
 put
\[ Q_j \defeq \frac{1}{\sqrt{N^n}}
\sum_{ \ell \in \ZZ^n } 
\delta( x - \ell - j/N) \; : \; j \in (\ZZ / N \ZZ)^n .Ä \]
We write the operator $ f^w ( x , hD ) $ on $ \Hh_h^n $ 
explicitly in that basis using the Fourier expansion of its symbol:
\[
 f ( x, \xi ) = \sum_{ \ell, m \in \ZZ^n } \hat f ( \ell, m ) \,
e^{ 2 \pi i ( \la \ell , x \ra + \la m , \xi \ra ) } \,.
\]
For that let $L_{\ell, m } ( x, \xi ) = \la \ell , x \ra + \la m , \xi \ra $,
so that
\[ 
f^w ( x, h D ) = \sum_{ \ell, m \in \ZZ^n } \hat f ( \ell, m ) \,
\exp( 2 \pi i L^w _{ \ell, m } ( x, h D ) ) \,. 
\]
We also check that 
\[
 \exp\big( 2 \pi i L^w _{ \ell, m } ( x, h D ) \big)\, Q_j = 
\exp ( {\pi i }( 2 \la j , \ell \ra - \la m , \ell \ra)/N )\,  Q_{j-m} \,,
\]
(note that $ j \in (\ZZ / N \ZZ)^n  $ and  $ j - m $ is meant $ \mod N $)
and consequently,
\begin{gather*}
 f_N ( Q_j )    = \sum_{ m \in \ZZ^n/( N \ZZ)^n } F_{mj} \, Q_m
\,, \\
F_{mj} =  \sum_{ \ell , r \in \ZZ^n  } 
\hat f ( \ell, j - m - rN) (-1)^{\la r, \ell \ra}\, 
\exp ( {\pi i}\, \la j + m , \ell \ra /N ) \,. 
\end{gather*}
Since
\[ 
\begin{split} 
\bar F_{jm} & = 
 \sum_{ \ell , r \in \ZZ^n  } 
\hat {\bar f} (- \ell, j - m + rN) (-1)^{\la r, \ell \ra} \,
\exp (- { \pi i }\,  \la j + m , \ell \ra /N ) \\ 
& = 
\sum_{ \ell , r \in \ZZ^n  } 
\hat {\bar f} ( \ell, j - m - rN) (-1)^{\la r, \ell \ra} \,
\exp ( { \pi i }\,  \la j + m , \ell \ra/N ) \,, 
\end{split} 
\]
we see that for real $ f $, $ f = \bar f $, $ F_{jm} = 
\bar F_{mj} $. This means that $ f^w ( x , hD ) $ is 
self-adjoint for the inner product $ \langle \bullet , \bullet \rangle_0 $. 
We also see that the map 
$ f \mapsto (F_{jm} )_{ j , m \in ( \ZZ / N \ZZ )^n } $ is onto, 
from $ \CI ( \tn ; \RR ) $ to the space of Hermitian matrices. 

Any other metric on $ \Hh_h^n $ could be written as $ \langle u , 
v \rangle  = \langle B u , v \rangle_0 = \langle u , B v \rangle_0 $.
If $ \langle f_N u , v \rangle = \langle u , f_N v \rangle $ for all real
$ f$'s, then $ B f_N= f_N B $ for all such $ f$'s, and hence for all
Hermitian matrices. That shows that $ B = c\, {\rm Id} $, as claimed.
\end{proof}

We normalize the inner product so that the basis specified in 
\eqref{eq:bas} is orthonormal. From now on we use this basis
to identify 
\[   \Hh_h^n \simeq \ell^2 ( \ZZ_N^n  ) \simeq \CC^{N^n} \,, 
\ \ \ZZ_N^n \defeq ( \ZZ / N \ZZ)^n \,. 
 \]

The calculation of the matrix coefficients in the 
proof of Lemma \ref{l:hilb} immediately gives the following
\begin{lem}
\label{l:trace}
Suppose $ f \in \CI ( \TT^{2n} ) $. Then 
\begin{equation}
\label{eq:trace2}
\begin{split}
\tr f_N & = N^n \sum_{ \ell , m \in \ZZ^n } ( -1)^{N\langle  \ell, m \rangle } \hat  f ( N \ell, N  m ) \\
& = N^n \int_{\TT^{2n} } f (w) d {\mathcal L} ( w)  + r_N \,,  \\
& \hspace{-0.5in} | r_N| \leq C_{kn} N^{-k+n } 
\sum_{ |\alpha | \leq \max(k, 2n+1) } \int_{\TT^{2n}} | \partial^{\alpha} f 
(w) |d {\mathcal L} ( w) 
 \,,
\end{split}
\end{equation}
for any $ k $. Here $ {\mathcal L }( w ) $ is the Lebesgue measure on 
$ \tn $ normalized so that $ { \mathcal L} ( \tn ) = 1 $.
\end{lem}

It is well known that for $ f \in \CI ( \tn ) $, independent of 
$ N $, $ f_N $ is uniformly bounded on $ \ell^2 ( \ZZ_N^n ) $ -- 
see \cite{BouzDB}. 
We will recall a slight generalization of that for functions
which are allowed to depend on $ N $ in a $ S_\rho $-way described
in \S \ref{pc}.

\subsection{$ S_\rho $ classes for the torus.}

The $ S_\rho $ classes for the quantization of the torus have 
already been considered in \cite{Sch} and we refer to that
paper for more detailed results such as the sharp G{\aa}rding
inequality. Here we continue with a self-contained presentation.

We first define a class of order functions: a function of 
$ w \in \TT^{2n} $ and $ \alpha > 0 $ 
is an $\alpha$-order function if there exist
$ C $ and $ M $ (independent of $ \alpha $) such that
\begin{equation}
\label{eq:epsor}
\forall \, w, w' \in \TT^{2n} \,, \ \ \
\frac{ m ( w, \alpha ) } { m ( w' , \alpha ) } \leq C 
( 1 + d_{ \TT_\alpha^n } ( w/\alpha , w'/\alpha ) )^M \,, \ \ 
\TT_\alpha^n \defeq (\RR^{2n}/(\ZZ/\alpha)^{2n}) \,, 
\end{equation}
with the distance induced from the Euclidean distance:
$ d_{ \RR^{2n} / \Gamma } ( w, w' ) = 
\inf_{ \gamma \in \Gamma} | w - w' + \gamma | $.

With this definition we have 
\begin{equation} 
\label{equation}
S ( m , \alpha ) \defeq \{ 
a \in \CI ( \TT^{2n} ) \,, \ \ \partial^\beta a ( w ) = {\mathcal O}
( \alpha^{-|\beta|} m ( w, \alpha ) ) \} \,. 
\end{equation}

If 
\begin{equation}
\label{eq:alN}  N^{-\rho} /C \leq  \alpha \leq C N^{-\rho} \,, \ \ 
0 < \rho < \frac 12 \,, \end{equation}
the quantization procedure described in \S \ref{qoc} applies to 
$ S ( m , \alpha ) $:
we now quantize functions $ f $ which are periodic and belong
to $ S_\rho$ with $ h = 1/(2 \pi N ) $. Similarly, we have the 
composition formula \eqref{eq:comptor} with the asymptotic
expansion in \eqref{eq:comp} valued in $ S ( m_1 m_2 , \alpha) $. 

Lemma \ref{l:ord} translates into this setting and will be used
in \S \ref{pot}:
\begin{lem}
\label{l:ord}
For $ f \in \CI ( \tn )   $ 
\[  m ( w, \alpha)  \defeq |f ( w )|^2 + \alpha^{2}  \,, \]
as an $ \alpha$-order function in the sense of definition \eqref{eq:epsor}.
In addition, 
for $  \psi \in \CIc ( \RR ; [ 0 , 1 ] ) $ equal to $ 1 $ on $ [ - 1, 1] $,
\begin{equation}
\label{eq:mell}   \left( |f ( w )|^2 + \alpha^2 \psi \left( 
\frac {| f ( w )|^2 }{ \alpha^2 } \right) \right)^{\pm 1} 
 \in S ( m^{\pm 1}  , \alpha ) \,. \end{equation}
\end{lem}

For $ S ( 1 , \alpha )$  we also have uniform $\ell^2 $-boundedness,
which we present in the simplest form:
\begin{prop}
\label{p:L2}
Suppose $ f \in S ( 1 , \alpha ) $ with $ \alpha $ satisfying 
\eqref{eq:alN}. Then 
\begin{equation}
\label{eq:l2b}     \| f_N \|_{ \ell^2 \rightarrow \ell^2 } \leq \sup_{\TT^{2n}} | f |  + 
o ( 1 )  \,, \ \ N \rightarrow \infty .
\end{equation}
\end{prop}
\begin{proof}
Lemma \ref{l:trace} gives
\[ \| f_N \|^2_{\rm{HS}} \defeq \tr f_N^* f_N = 
N^n \int_{\tn} \bar f \#_h f d \cL + {\mathcal O} ( N^{-k+n} ) 
\sum_{ |\beta| \leq k  } \int_{\tn } | \partial^\beta ( \bar f \#_h f) | 
d \cL \,, \ \ k \gg n \,.  \]
Since $ \bar f \#_h f \in S( 1 , \alpha ) $ (that is, using \eqref{eq:alN},
$ \bar f \#_h f$
lies in $ S_\rho $ when considered
 as a periodic
function on $ \RR^{2n} $),
we see that
\[ \| f_N \|^2_{\rm{HS}} = {\mathcal O} ( N^n ) + 
{\mathcal O}  ( N^{-k( 1 - \rho ) + n} ) = {\mathcal O} ( N^n) \,. \]
Hence, 
\begin{equation}
\label{eq:HSf}
 \| f_N \|_{ \ell^2 \rightarrow \ell^2 } \leq \| f_N\|_{\rm{HS}} \leq C N^{\frac n2} 
\,.
\end{equation}
We now use H\"ormander's trick for deriving $L^2$-boundedness from 
the semiclassical calculus. 
Let $ M > \sup_{\tn}  | \bar f \#_h f | $ and let $ 
a_N = M - f_N^* f_N $, $ a = M - \bar f \#_h f \in S ( 1, \alpha ) $, $ a > 
1/C > 0 $. 
Then by \eqref{eq:comp}
\[   b^0_N b^0_N - a_N = r^0_N \,, \ \ r^0 \in N^{2 \rho-1} S ( 1 , \alpha) \,,
\ \  b^0 \defeq  \sqrt a \in S ( 1 , \alpha ) \,. \]
We now proceed by induction to construct real $ b^j \in N^{ j ( 2 \rho - 1 ) }  
S ( 1, \alpha ) 
 $, $ 0 < j \leq J $,
so that 
\begin{equation}
\label{eq:Hind}  (B^J_N)^2  - a_N = r^J_N \,, \ \ 
B^J_N \defeq \sum_{j=0}^{J} 
b^j_N \,, \ \ 
r^J \in N^{(J+1)(2\rho-1)} S ( 1 , \alpha )  \,. 
\end{equation}
Suppose that we already have it for $ J $ (the first inductive 
step being $ J=0$) and we want to find $ b^{J+1} \in  N^{ ( J + 1 ) 
( 2 \rho - 1 ) }  
S ( 1, \alpha ) $ so that : 
\[ \begin{split} \left( B_N^J +  
b^{J+1}_N \right)^2 - a_N & = r_N^J +  B_N^J b_N^{J+1} + b_N^{J+1} B_N^J 
\\ & 
= r_N^J + b_N^0 b_N^{J+1} +  b_N^{J+1} b_N^0 + R_N^J b_N^{J+1} +
b_N^{J+1} R_N^J \,, 
\end{split} 
\]
where $ R^J = B^J - b_0 \in N^{2\rho-1} S( 1 , \alpha ) $. We now
simply put 
$$ b^{J+1} = r^J/(2 b^0 ) \in  N^{(J+1)(2\rho-1)} S ( 1 , \alpha )  \,,
$$
which is real since the left hand side of \eqref{eq:Hind} is
self-adjoint. The inductive step follows again from the composition 
property.

Returning to the boundedness on $ \ell^2 $ we now have
\[ \begin{split}
M \| u \|^2 - \| f_N u \|^2 &  = \langle a_N u , u \rangle 
 = \langle B_N^J u , B_N^J u \rangle - \langle r_N^J u , u \rangle 
\\
& \geq - \| r_N^J \|_{\ell^2 \rightarrow \ell^2 }  \| u \|^2 \geq -
 \| r_N^J \|_{\rm{HS}} \| u \|^2 \\
& \geq - C N^{\frac n2 + (J+1)(2 \rho -1) } \| u \|^2 \,,
\end{split}
\]
where for the last inequality we used \eqref{eq:HSf}. Hence by
taking $ J $ large enough, 
$ \| f_N \|_{\ell^2 \rightarrow \ell^2} \leq M^{1/2} + o(1) $, and
since $ M $ can be taken as close to $ \sup |f| $ as we like,
this gives \eqref{eq:l2b}
\end{proof} 

One of the consequences of the boundedness on  $\ell^2 $ is 
the justification of the basic principle of semiclassical 
quantization: 
\[ \text{ Poisson brackets, $\{ \bullet, \bullet \}$ } \ \longleftrightarrow 
\text{ Commutators, $ (i/h) [ \bullet_N , \bullet_N ] $, \ \ $ h = 1/(2 \pi N ) $.}\]
More precisely, $ \{ f , g \} = \sum_{ j=1}^n ( \partial_{\xi_j} f 
\partial_{x_j} g - \partial_{x_j} f  \partial_{\xi_j} g) $, (with
$ \sigma = \sum_{j=1}^n d\xi_j \wedge dx_j $ the symplectic form 
on $ \tn $), and 
\begin{equation}
\label{eq:po2co}
2 \pi i N [ f_N , g_N] = \left( \{ f , g \} \right)_N + {\mathcal O}_{\ell^2
\rightarrow \ell^2 }  ( N^{-2+4\rho} ) \,. 
\end{equation}

The functional calculus lemma presented in the $ \RR^n $ setting
translates to the case of the torus:
\begin{lem}
\label{l:func}
Suppose $ f \in \CI ( \TT^{2n} ) $ and $ \alpha = h^\rho $,  $0 \leq \rho < 1\
 /2 $.
Then, for $ \psi \in \CIc ( \RR ) $,
\begin{gather}
\label{eq:func}
\begin{gathered}
\psi \left( \frac{f_N^* f_N } {\alpha^2 } \right) = q_N \,, \ \
q \in S ( 1 , \alpha )\,, \\
q = q^0 + h^{1- 2 \rho} q^1 + {\mathcal O}_{S } ( h^\infty) 
\,, \ \ q^j \in  S ( 1 , \alpha )\,, \ \
q^0 ( w )  = \psi ( |f ( w ) |^2/ \alpha^2 ) \,, \\
q^1 ( w ) = \widetilde \psi ( |f ( w ) |^2/ \alpha^2 ) 
\widetilde q^1 ( w) \,, \ \ \widetilde q^1 \in S( 1 , \alpha ) \,, \ \ 
\widetilde \psi \in \CIc ( \RR ) \,, \ \ \widetilde \psi\rest_{{\rm{supp}}\, \psi} 
\equiv 1 \,. 
\end{gathered}
\end{gather}
\end{lem}
\begin{proof}
We need to check that for a function $ \varphi \in \CIc ( \RR ) $, and
$ g \in \CI ( \tn; \RR ) $, the action of $ \varphi ( g_N ) $ on 
$ \Hh_h^n \simeq \ell^2 ( \ZZ_N^n ) $ defined using functional
calculus of self-adjoint matrices is the same as the action of
$ \varphi( g^w ( x, h D ) ) $ on $ \Hh_h^n \subset \cS' ( \RR^n ) $.
In view of the Helffer-Sj\"ostrand formula that follows 
from verifying that the action of the resolvent $ ( z - g_N )^{-1} $,
$ \Im z \neq 0 $, on $ \Hh_h^n $ is the same as the 
action of $ ( z - g^w ( x , h D ) )^{-1} $, $ \Im z \neq 0 $,
on $ \Hh_h^n $ as a subset of $ \cS ( \RR^n ) $. But we know 
from \eqref{eq:Beals} that for $ \Im z \neq 0 $,
$ ( z - g^w ( x , h D ) )^{-1} = F(z, x, h D ) $, where $ F ( z) \in 
S ( 1 ) $ (non-uniformly as $ \Im z \rightarrow 0 $ but with 
seminorms polynomially bounded). This means that the $ L^2 $ inverse
is a restriction of an inverse defined on $ \cS' ( \RR^n ) $.
Hence $ ( z - g_N )^{-1} = [F( z)]_N $ and the actions are the same.
This argument is not asymptotic in $ N $ and 
applies to $ \varphi = \psi ( \bullet / \alpha^2 ) $ and $ g = \bar f \#_h f$.
\end{proof}

\begin{prop}
\label{p:rank}
Suppose that \eqref{eq:cond} holds with $ z = 0 $. Then 
for any $ \psi \in \CIc ( \RR ) $,
\begin{equation}
\label{eq:rank}
\rank \psi \left ( \frac{ f_N f_N^*} {\alpha^2} \right)  \leq 
C N^{n } \alpha^{ \kappa } 
\,, \ \ N^{-\rho} \leq \alpha \ll 1 \,, \ \
\rho < \frac 12  \,,
\end{equation} 
with the constant depending only on the support of $\psi $.
\end{prop}
We note that by 
proceeding 
either as in the proof of \cite[Proposition 4.4]{h-sj}
or as in the proof of \cite[Proposition 5.10]{SjZw04}
we can show that the result is valid for $ \rho = 1/2 $ but 
we do not need that in this paper.
\begin{proof}
Suppose $ \psi_1 \in \CIc ( ( - R^2 +1, R^2-1 ) ,[0,1] ) $, $ R \gg 1 $ 
is equal to $ 1 $ on 
the support of $ \psi $. Then, using the functional calculus of self-adjoint
matrices and Lemmas \ref{l:trace}, \ref{l:func}, and \eqref{eq:comp} 
we get, with $ \widetilde \psi \in \CIc ((-R^2,R^2),[0,1])$, $ \widetilde \psi\rest_{ {\rm{supp}
} \, \psi_1 } \equiv 1 $, 
\[  \begin{split}
\rank \psi \left ( \frac{ f_N f_N^*} {\alpha^2} \right)  & \leq
\tr \psi_1 \left ( \frac{ f_N f_N^*} {\alpha^2} \right)  
 \leq N^n \int_{\tn} \widetilde \psi ( |f|^2 / \alpha^2 ) 
 d {\mathcal L } + {\mathcal O} (N^{-\infty } ) \\
& \leq N^n {\mathcal L} ( \{ w :  |f ( w ) | \leq R \alpha \} ) 
+  {\mathcal O} (N^{-\infty } )
 \leq C  N^{n } \alpha^{ \kappa } \,, 
\end{split}
\]
proving the lemma.
\end{proof}

\section{Some facts about random matrices}
\label{rmt}

\def \rank {\operatorname{rank}}
\def \comp {\operatorname{comp}}
\def \supp {\operatorname{supp}}
\def \diam {\operatorname{diam}}
\def \dist {\operatorname{dist}}
\def \tr {\operatorname{tr}}
\def \Vol {\operatorname{Vol}}
\def \mco {{\mathcal O}}
\def \Domain {\operatorname{Domain}}
\def \vol {{\rm vol}}
\def \pilambda {\Pi^{\partial X_0}_{\Lambda}}
\def \pilambdanot {\Pi^{\partial X_0}_{\Lambda_0}}
\def \piN {\Pi^{\rm{in}}_{N}}

\def \rest {\!\!\upharpoonright}
\def \restrict {\!\!\upharpoonright}
\def \mfm {{\mathfrak M}}
\def \mfp {{\mathfrak P}}
\def \mfg {{\mathfrak G}}
\def \mfo {{\mathfrak O}}
\def \mcd {{\mathcal D}}
\def \Real {{\mathbb R}}
\def \RR {{\mathbb R}}
\def \CC {{\mathbb C}}
\def \Sphere {\mathbb{S}}
\def \Complex {\mathbb{C}}
\def \Natural {{\mathbb N}}
\def \Sphere {{\mathbb S}}
\def \complement {{\bf C}}
\def \mcr{{\mathcal R}}
\def \hol {{\mathcal H}}
\def \Integers {{\mathbb Z}}

\def \mcl {{\mathcal L}}
\def \mcalp {{\mathcal P}}

\def \pr {{\mathcal J}}

\def \Ev {{\mathbb E}}
\def \Pr {{\mathbb P}}
\def \bbone {{\mathbb 1}}

\renewcommand\thefootnote{\dag}%

Random matrix theory is a very active field and we refer to
Mehta's classic book \cite{Me} for general background, 
and to \cite{FSV} for some
recent works and applications.
All the facts we need in this paper are elementary but they do not seem 
directly present in the mainstream literature. 
Consequently the presentation is almost self-contained
 and, reflecting
the authors's own position, does not assume any knowledge of the subject.

We consider the ensemble of complex Gaussian matrices with independent
entries distributed in $ \CC $ according to the standard normal
distribution. That means that there exists a probability space,
$ (\Omega , \Sigma , \mu ) $, $ \Sigma $ a $ \sigma$-algebra of
subsets of $ \Omega $ and $ \mu : \Sigma \rightarrow [0, \infty ) $,
a measure, with $ \mu ( \Omega ) = 1 $, and
a map $ \Omega \ni \omega \mapsto 
A_N ( \omega ) $, $ A_N ( \omega )  = ( a_{ij} ( \omega ) )_{ 1 \leq
i, j \leq d} $, such that $ \omega \mapsto a_{ij} ( \omega ) $ 
are {\em independent} random variables with standard normal distribution.
The push forward measures on $ \CC $, $ (a_{ij})_* \mu $,
are given by $ \exp(-|z|^2) d {\mathcal L} ( z ) / \pi $, where
$ {\mathcal L} $ is the Lebesgue measure (standard normal distribution),
and 
\[  [(a_{ij},a_{k\ell})]_* \mu = 
\frac{ 1}{  \pi^2 } e^{  -  |z|^2 - |w|^2 }  d {\mathcal L} ( z )
d {\mathcal L} ( w) \,, \ \ (i,j) \neq ( k , \ell) \,, \]
$ ( a_{ij} , a_{k \ell} ) : \Omega \rightarrow \CC_z \times \CC_w $,
which is the statement that $ a_{ij} $ and $ a_{k \ell} $ are independent.

A more useful global description of the random variable $ A_d ( \omega ) $
is given as follows: let $a_i= (a_{i1},...,a_{id})^t\in \CC^d$,
and set $A=(a_1,...,a_d)$.
Denote 
 $$d\mcl (a_{i})= d \Re a_{i1} d\Im a_{i1}...d\Re a_{id} d\Im a_{id}, 
\text{ and }  
d\mcl(A)= \prod_{i=1}^d d\mcl( a_{i}).$$ Then, as a measure on $ \CC^{d^2} $,
the space of $ d \times d $ matrices, 
\begin{equation}
\label{eq:LA}
A_* \mu  = {   \pi ^{-d^2} } \exp \left( - 
\| A \|_{\rm{HS}}^2 \right)d{\mathcal L}(A) \,, \ \ \| A \|^2_{\rm{HS}}  \defeq 
\tr A^* A \,,
\end{equation}
where HS stands for Hilbert-Schmidt.  Note that each entry $a_{ij}$ of $A$ 
is a complex $N(0,1)$ random variable.

We recall that any matrix $ A $ can be written using its singular value
decomposition, 
\begin{equation}
\label{eq:svd}  A = U S V^* \,, 
\end{equation}
where $ U U^* = U^* U = Id $, $ V V^* = V^* V = Id $, that is $ U $ and 
$ V $ are unitary, and $ S $ is a diagonal matrix with non-negative
entries. If the entries of $ S $ are distinct and we order them,
the decomposition is unique.

\begin{figure}
\begin{center}
\includegraphics[width=6in]{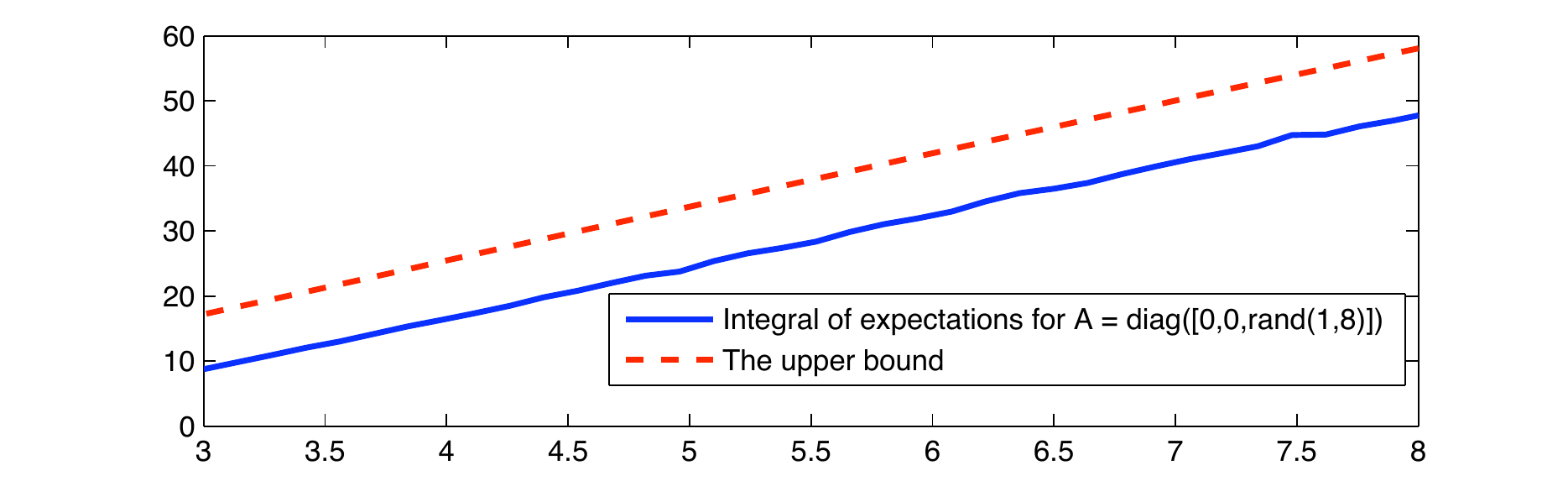}
\end{center}
\caption{A numerical example suggesting that Proposition \ref{p:tA} is
optimal: the left hand side is computed numerically for {\tt
A = diag([0,0,rand(1,8)]) } (a $10\times10$ diagonal matrix of rank $ 8 $)
 where {\tt rand} command produces
uniform distribution on $ [0,1]$. It is plotted
as a function of $ \log(1/\delta ) $. The upper bound in
Proposition \ref{p:tA} (with $ C = 1$) is also plotted for
comparison.}
\label{f:4}
\end{figure}

\begin{prop}\label{p:tA}
Let $A$ be a constant $d\times d$ matrix, and 
let  $Q$ be  a $d\times d$ random
matrix, with the entries $q_{ij}$ independent 
complex $N(0,1)$ random variables. Then there exists a constant $ C $
independent of $d$ and 
$A$, such that
 $$   \int_0^1 \left| \Ev (\tr(tA+\delta Q)^{-1}A) \right| dt 
\leq 
C \tr \left( \frac {|A| }{ \delta + |A| } \log \left( 2 + \frac{|A|} \delta 
\right) \right) \,, 
$$
where $ | A | = \sqrt { A A^* } $.
\end{prop}
In the proof we will need the following
\begin{lem}\label{l:g}
The function $g(s) \defeq \int_{\CC}  (1/|s+q|) e^{-|q|^2}d{\mathcal L}(q)$ is
smooth for $s\in \Complex$, and 
\[  g ( s ) = \frac{ \pi } { | s| } + {\mathcal O} \left( \frac 1 { | s|^2} 
\right) \,, \text{as } |s| \rightarrow \infty. \]
\end{lem}
\begin{proof}
The asymptotic expansion follows from the integrability of $ | q | $,
a change of variables, $ q \mapsto 1 + q/s $, and the method of 
stationary phase. 
\end{proof}

\begin{proof}[Proof of Proposition \ref{p:tA}]
Using the singular value decomposition for $A$, we may write $A=U S V^*$,
with $U$, $V$ unitary  and $S$ a diagonal matrix with non-negative
entries $\sigma_1,...,\; \sigma_d$ on the diagonal.  We note that 
$$\tr( (tA+\delta Q)^{-1}A)=\tr((t U S V^* +\delta Q)^{-1} U S V^* )
=  \tr ( (tS+\delta U^* Q V)^{-1} S).$$ 
Since $U^* Q V $ is a random matrix with the same probability 
distribution function as $Q$,
we have 
$$\Ev (\tr((tA+\delta Q)^{-1}A))= \Ev (\tr (( t S +\delta Q)^{-1} S )).$$
Thus we may assume that $A$ is diagonal, with non-negative entries
$\sigma_1,\; ...,\; \sigma_d$.
\renewcommand\thefootnote{\ddag}%
We have
$$\tr( (tA+\delta Q)^{-1}A)= \sum_1^d 
\frac {M_{ii}\sigma_i}{\det (tA+\delta Q)} $$
where here and below $M_{ij}$ is the $(i,j)$ minor of the matrix 
$tA+\delta Q$.  

To compute
$\Ev  ( {M_{ii}\sigma_i}/{\det (tA+\delta Q)})  $, 
we write 
$$\det(tA+\delta Q) =(t\sigma_i +\delta q_{ii})M_{ii}
+ \sum_{j\not = i} (-1)^{j+i} \delta q_{ij} M_{ij}$$
and define 
\begin{equation}
\label{eq:Sigmaii}
\Sigma_{ii} \defeq \left\{ q \in \CC^{d^2} \; : \; 
|(t\sigma_i+\delta q_{ii})M_{ii}| 
> \left|\sum_{j\not = i} (-1)^{j+i} \delta q_{ij} M_{ij}\right| \right\} \,. 
\end{equation}
Let $\bbbone_{F}$ be the characteristic function of a set $F$.
Then
\begin{equation}\label{eq:sepexp}
\Ev  \left(\frac {M_{ii}\sigma_i}{\det (tA+\delta Q)}\right)
= \Ev  \left(\frac {M_{ii}\sigma_i}{\det (tA+\delta Q)}\bbbone_{\Sigma_{ii}}\right)
+\Ev  \left(\frac {M_{ii}\sigma_i}{\det (tA+\delta Q)}
\bbbone_{\overline{\Sigma}_{ii}^{c}} \right)
\end{equation}
since the boundary of $\Sigma_{ii}$ has measure $0$\footnote{This follows
from the fact that the pushforward of the probability measure 
by $ Q $ (the probability density) is absolutely continuous with 
respect to the Lebesgue measure on $ \CC^{n^2} $ and the set 
$$ 
\{ Q \in \CC^{n^2} \; : \;  Q = ( q_{ij} )_{1\leq i,j \leq n} \,, \ \
|(t\sigma_i+\delta q_{ii})M_{ii}|^2
= |\sum_{j\not = i}^d (-1)^{j+i} \delta q_{ij} M_{ij}|^2 \} \,, $$
has Lebesgue measure $ 0 $.}.

Now,
\[ \begin{split}
\Ev  \left(\frac {M_{ii}\sigma_i}{\det (tA+\delta Q)}\bbbone_{\Sigma_{ii}}\right)
&
= 
\Ev  \left(\frac {M_{ii}\sigma_i}{ (t\sigma_i +\delta q_{ii})M_{ii}}
\left( 1+  \frac{ \sum_{j\not = i}^d (-1)^{j+i} \delta q_{ij} M_{ij}}
{(t\sigma_i +\delta q_{ii})M_{ii}} \right)^{-1} \bbbone_{\Sigma_{ii}}\right)\\
&
= \Ev  \left(\frac {\sigma_i}{ (t\sigma_i +\delta q_{ii})}
\sum _{k=0}^{\infty}\left( - 
\frac{ \sum_{j\not = i}^d (-1)^{j+i} \delta q_{ij} M_{ij}}{
           (t\sigma_i +\delta q_{ii})M_{ii}} \right)^{k}\bbbone_{\Sigma_{ii}}\right).
\end{split} \]
We recall that the set $\Sigma_{ii}$ is chosen so that the infinite 
sum converges. 

The set $\Sigma_{ii}$ is invariant under the mapping
\begin{equation}
\label{eq:qimapping}
q_{i1},...,q_{i,i-1}, q_{i,i+1},...,q_{i,d}\mapsto
e^{i\varphi} q_{i1},...,e^{i\varphi} q_{i,i-1}, e^{i\varphi}  q_{i,i+1},...,
e^{i\varphi} q_{i,d} 
\end{equation}
for any real number $\varphi$.  Since $ M_{ij} $'s are independent
of $ q_{ij} $, 
$  \sum_{j\not = i}^d (-1)^{j+i} \delta q_{ij} M_{ij} $
is homogeneous of degree $1$ under this same mapping and 
$(t\sigma_i +\delta q_{ii})M_{ii}$ is independent of $q_{ij}$ for $j\not = i$,
 we find that 
$$\Ev  \left(\frac {M_{ii}\sigma_i}{\det (tA+\delta Q)}\bbbone_{\Sigma_{ii}}
\right)
= \Ev  \left(\frac {\sigma_i}{ (t\sigma_i +\delta q_{ii})}
\bbbone_{\Sigma_{ii}}\right).$$

We do a similar computation for the second term of (\ref{eq:sepexp}):
\begin{gather*}
  \Ev  \left(\frac {M_{ii}\sigma_i}{\det (tA+\delta Q)}
\bbbone_{\overline{\Sigma}_{ii}^{c}} \right)
= \Ev  \left(\frac {M_{ii}\sigma_i}
{\sum_{j\not = i}^d (-1)^{j+i} \delta q_{ij} M_{ij}} 
\left( 1+ \frac{ (t\sigma_i+\delta q_{ii}) }
{\sum_{j\not = i}^d (-1)^{j+i} \delta q_{ij} M_{ij}}\right)^{-1}
\bbbone_{\overline{\Sigma}_{ii}^{c}} \right)
\\ 
= \Ev  \left(\frac {M_{ii}\sigma_i}
{\sum_{j\not = i}^d (-1)^{j+i} \delta q_{ij} M_{ij}} \sum_{k=0}^{\infty}
\left( - \frac{ (t\sigma_i+\delta q_{ii}) }
{\sum_{j\not = i}^d (-1)^{j+i} \delta q_{ij} M_{ij}}\right)^{k}
\bbbone_{\overline{\Sigma}_{ii}^{c}} \right)
= 0\,,
\end{gather*}
using, as before, the invariance properties of $\Sigma_{ii}$ and the 
homogeneity of 
$$ {\sum_{j\not = i}^d (-1)^{j+i} \delta q_{ij} M_{ij}} \,. $$

Thus we have 
\begin{equation}
\label{eq:traceA}
 \Ev (\tr(tA+\delta Q)^{-1}A) 
= \sum _{i=1}^d \Ev  \left(\frac {\sigma_i}{ (t\sigma_i +\delta q_{ii})}
\bbbone_{\Sigma_{ii}}\right).
\end{equation}

Now,
\[ \begin{split}
\left|\int_0^1 \Ev  \left( \frac {\sigma_i}{ (t\sigma_i +\delta q_{ii})}
\bbbone_{\Sigma_{ii}} \right) dt \right|
& 
\leq \int _0^1 
\Ev \left( \frac{\sigma_i}{|t \sigma_i +\delta q_{ii}|}\right) dt  
= \int_0^1\Ev \left( \frac{\sigma_i/\delta}{|t \sigma_i/\delta +q_{ii}|}
\right) dt
\\ & = \int_0^{\sigma_i/\delta}\Ev \left(\frac{1}{|s+q _{ii}|}\right) ds 
 = \frac{1}{\pi} \int_0^{\sigma_i/\delta}g(s)ds
\end{split} \] 
where $g$ is the function defined in Lemma \ref{l:g}.  Using this,
(\ref{eq:traceA}), and the results of Lemma \ref{l:g} proves the 
proposition. 
\end{proof}

\begin{lem}\label{l:invertible}
Let $F,\; G$ be $d\times d$ matrices, with $F$ invertible, and let $\beta 
= \|F^{-1}\|.$
Then
$$
\Ev  \Big(\tr \big( (F+\delta Q)^{-1}G \big)\Big) 
=  \tr \Big( F^{-1} G \Big)\big(1+ O(e^{-1/4(\delta \beta d)^2})\big)+ 
{\mathcal O}\left( \frac{1}{\delta}
\|G\|d^4 e^{-1/4 (d\beta \delta)^2}\right).
$$
The implicit constant in the error term is independent of 
$F$ and $G$.
\end{lem}
\begin{proof}
We first note that if we replace $F$ by its singular value decomposition,
$F=U S V^*$, then
$$\Ev  \Big(\tr \big( (F+\delta Q)^{-1}G \big)\Big) 
= \Ev  \Big(\tr \big( (S +\delta Q)^{-1}(U^*GV) \big)\Big)$$
and 
$$\tr \big( F^{-1} G \big)= \tr \big(S^{-1} U^*G V \big).$$
Thus we may assume that $F$ is a diagonal matrix.

Our proof then resembles the proof of Proposition \ref{p:tA}.
Let $\chi\in L^{\infty}(\Real_+)$ be the characteristic function of 
$(- \infty , 1/2 ] $,
and, if $A=(a_{ij})$, let $\| A\|_{\sup}= \sup_{ij} |a_{ij}|.$
We write
\begin{multline}\label{eq:splitev}
\Ev\Big(\tr \big( (F+\delta Q)^{-1}G\big) \Big) = 
\Ev\Big(\tr \big( (F+\delta Q)^{-1}G\big)\chi(d \|Q\|_{\sup}\delta \beta)\Big) \\
+ \Ev\Big(\tr \big((F+\delta Q)^{-1}G\big)
(1-\chi(d \|Q\|_{\sup}\delta \beta))\Big)
.
\end{multline}
For the first term,
\begin{align*}
\Ev\Big(\tr\big( (F+\delta Q)^{-1}G\big) \chi(d \|Q\|_{\sup}\delta \beta)\Big)
= \Ev \Big(\tr \Big (F^{-1} \sum_0^{\infty} 
(-\delta (Q F^{-1} )^j G \Big) \chi(d \|Q\|_{\sup}\delta \beta))\Big).
\end{align*}
Using the fact that the cut-off $\chi(d \|Q\|_{\sup}\delta \beta)$ is 
invariant under rotations of the $q_{ij}$ and that the $q_{ij}$ are complex
and independent, we find
\begin{align}\label{eq:part1} 
\Ev\Big(\tr \big( (F +\delta Q)^{-1} G \big) 
\chi(d \|Q\|_{\sup}\delta \beta)\Big)
& = \tr\big ( F^{-1} G ) \mu ( Q :\|Q\|_{\sup}< 1/2\delta \beta d) \nonumber 
\\
& = \tr\big ( F^{-1}G ) (1+ {\mathcal O}(d^2e^{-1/4(\delta \beta d)^2})) \,. 
\end{align}

Now we consider the remaining term of (\ref{eq:splitev}).
In a way similar to the proof  of Proposition \ref{p:tA}, 
we denote the diagonal entries of $F$ by $f_{ii}=\sigma_i$, and 
by $M_{ij}$ the $(i,j)$ minor of $F+\delta Q$.  If $G=(g_{ij})$, we 
have 
$$
\Ev\Big(\tr \big((F+\delta Q)^{-1}G\big)
(1-\chi(d \|Q\|_{\sup}\delta \beta))\Big)
= \Ev \Big(\sum_{i,j} \frac{(-1)^{i+j}M_{ji}g_{ji}}{\det(F+\delta Q)}
(1-\chi(d \|Q\|_{\sup}\delta \beta))\Big).
$$
Just as in the proof of Proposition \ref{p:tA}, to compute
$$\Ev \Big( \frac{M_{ii}g_{ii}}{\det(F+\delta Q)}
(1-\chi(d \|Q\|_{\sup}\delta \beta))\Big)
$$ we write 
$$\det(F+\delta Q) =(\sigma_i +\delta q_{ii})M_{ii}
+ \sum_{j\not = i} (-1)^{j+i} \delta q_{ij} M_{ij}$$
and define $\Sigma_{ii}$ as in (\ref{eq:Sigmaii}).
Proceeding almost exactly as in the proof of Proposition \ref{p:tA}, using 
that both $\Sigma_{ii}$ and 
the support of $(1-\chi(d \|Q\|_{\sup}\delta \beta))$ are 
 invariant under the mapping (\ref{eq:qimapping}),
we get that 
\begin{equation*}
 \Ev \Big( \frac{M_{ii}g_{ii}}{\det(F+\delta Q)}
(1-\chi(d \|Q\|_{\sup}\delta \beta))\Big)
= \Ev \left( \frac{ g_{ii}}{(\sigma_ i + \delta q_{ii})}
\bbbone_{\Sigma_{ii}} (1-\chi( d \|Q\|_{\sup}\delta \beta))\right) .
\end{equation*}
But 
\begin{equation*}
\left|\Ev \left( \frac{ g_{ii}}{(\sigma_i + \delta q_{ii})}
\bbbone_{\Sigma_{ii}} (1-\chi( d \|Q\|_{\sup}\delta \beta))\right)\right|
\leq 
C \frac{\|G\|}{\delta} d^2 e^{-1/4(d\delta \beta)^2}.
\end{equation*}

To compute 
$$\Ev \Big( \frac{(-1)^{i+j}M_{ji}g_{ji}}{\det(F+\delta Q)}
(1-\chi(d \|Q\|_{\sup}\delta \beta))\Big)$$
when $i\not = j$,
we write
\begin{equation}\label{eq:detexp}
\det(F+\delta Q) =\delta q_{ji}M_{ji}(-1)^{i+j}
+ (\sigma_i+\delta q_{ii})M_{ii}
+ \sum_{k\not = i,j} (-1)^{k+i} \delta q_{ki} M_{ki}
\end{equation}
and define
$$\Sigma_{ji}\defeq \left\{ q \in \CC^{d^2} \; : \; 
|\delta q_{ji}M_{ji} | 
> \left|(\sigma_i+\delta q_{ii})M_{ii}+
\sum_{ i \neq k \neq j} (-1)^{k+i} \delta q_{ki} M_{ki}\right| \right\} \,.
$$
Following the proof of Proposition \ref{p:tA} but
treating the term $\delta q_{ji}M_{ji}$ as the 
distinguished one in the expansion of the determinant (\ref{eq:detexp}) and  using the invariance
of $\Sigma_{ji}$ under rotations of $q_{ji}$, we find that
\begin{multline*}
\Ev \left( \frac{(-1)^{i+j}M_{ji}g_{ji}}{\det(F+\delta Q)}
(1-\chi(d \|Q\|_{\sup}\delta \beta))\right) \\
= \Ev \left( \frac{(-1)^{i+j}M_{ji}g_{ji}}
{ (\sigma_i+\delta q_{ii})M_{ii}
+ \sum_{k\not = i,j} (-1)^{k+i} \delta q_{ki} M_{ki}}
\bbbone_{\overline{\Sigma}_{ji}^{c}}(1-\chi(d \|Q\|_{\sup}\delta \beta))\right).
\end{multline*}
Since on 
the support of $\bbbone_{\overline{\Sigma}_{ji}^{c}}$
$$|M_{ji}| \leq \frac{1}{\delta |q_{ji}|} \left| (\sigma_i+\delta q_{ii})M_{ii}
+ \sum_{k\not = i,j} (-1)^{k+i} \delta q_{ki} M_{ki})\right|
$$ 
we find 
$$
\left| \Ev \left( \frac{(-1)^{i+j}M_{ji}g_{ji}}{\det(F+\delta Q)}
(1-\chi(d \|Q\|_{\sup}\delta \beta))\right)\right|
\leq C \frac{\|G\| }{\delta} d^2 e^{-1/4(d\delta \beta)^2}.$$

\end{proof}

Our proof of Proposition \ref{p:ra2de} in the next section
will use Proposition \ref{p:trick}.
To prove this proposition we will need several preliminary results.

The first lemma below follows from well-known 
 facts about eigenvalues of complex Gaussian ensemble. We
give a direct  proof suggested to us by Mark Rudelson:
\begin{lem}\label{l:L1}
Let $A=(a_1,...,a_d)$, with $a_i\in \Complex^d$. Then,
with the notation of \eqref{eq:LA},
$$\int_{\| A\|_{HS}\leq 1} |\det  A|^{-1} d\mcl(A) <\infty\,. $$
\end{lem}
\begin{proof}
We begin by introducing some more notation.  
For $p \leq d$, $p\in \Natural$, $v\in \Complex^d$
denote by $\mcalp_p v$ projection onto the 
subspace spanned (over the complex numbers) by $a_1,...,a_p$.  This of 
course depends on $a_1,...,a_p$, but we omit this in our notation for 
simplicity.

Using the Graham-Schmidt process, we can, if $A$ is invertible (as 
it is off a set of measure $0$),
write the matrix $A=UR$, with $U$ a unitary matrix and $R$ being 
upper triangular.  The diagonal entries of $R$ are then given by 
$\|a_1\|$ and 
$\|(1 -\mcalp_{p-1})a_p\|$, $p=2,...,d$.  Thus 
\begin{align*}
|\det A|= \|a_1\|\| (1-\mcalp_{1})a_2\|
  \| (1-\mcalp_{2})a_3\| \cdot \cdot \cdot      \| (1-\mcalp_{d-1})a_d\|.
\end{align*}
Note that 
$$\|a_1\|\| (1-\mcalp_{1})a_2\|
  \| (1-\mcalp_{2})a_3\| \cdot \cdot \cdot \| (1-\mcalp_{d-2})a_{d-1}\|
$$
is independent of $a_d$, that is, independent of $a_{1d},a_{2d},...,a_{dd}$. 
Therefore
\begin{align*} & 
\int_{\| A\|_{HS}\leq 1} |\det  A|^{-1} d\mcl(A)
\\ & 
= \int _{\| A\|_{HS}\leq 1}
\frac{1}
{ \|a_1\|\| (1-\mcalp_{1})a_2\|
 \cdot \cdot \cdot \| (1-\mcalp_{d-2})a_{d-1}\|}
\frac{d{\mathcal L}( a_d) d{\mathcal L}( a_{d-1})...d{\mathcal L}( a_{1}) } 
{ \| (1-\mcalp_{d-1})a_d\|} \\ &
\leq \int_{\|a_1\| \leq 1} \cdot \cdot \cdot 
\int_{\|a_d\| \leq 1} \frac{d{\mathcal L}( a_d) } { \| (1-\mcalp_{d-1})a_d\|} 
\frac{ d{\mathcal L}( a_{d-1})...d{\mathcal L}( a_{1})}
{ \|a_1\|\| (1-\mcalp_{1})a_2\|
 \cdot \cdot \cdot \| (1-\mcalp_{d-2})a_{d-1}\|}.
\end{align*}
The value of 
$\int_{\|a_d\| \leq 1} 
1/{\| (1-\mcalp_{d-1})a_d\|}
d{\mathcal L}( a_d)$
depends only on $d$ and the rank of the space spanned by $a_1,...,a_{d-1}$.
We find
 $1/\| (1-\mcalp_{d-1})a_d\|$ is locally integrable over 
$\Real^{2d}\simeq \Complex ^d$, because 
$a_d\in \Complex^d$ and the space spanned by $a_1,...,a_{d-1}$ has complex
dimension at most $d-1$.Therefore
\begin{equation}\label{eq:ad}
  \int_{\|a_d\| \leq 1} \frac{1}{\| (1-\mcalp_{d-1})a_d\|}
d{\mathcal L}( a_d) \leq C <\infty.
\end{equation}
Here the constant $C$ can be chosen independent of 
$a_1,...,a_{d-1}$, as the maximum of the integral in (\ref{eq:ad}) occurs 
when $a_1,...,a_{d-1}$ span a $d-1$ dimensional vector space.
The proof follows by iterating the above argument.
\end{proof}

\begin{prop}\label{p:trick}
Let $A(s,t)$ be a $d\times d$ matrix depending smoothly on 
$(s,t)\in U \subset \Complex^2$. 
Let $Q$ denote a $d\times d$ random matrix, with each entry 
an independent complex $N(0,1)$ random variable. Then for $\delta >0,$
$(s,t)\in U$,
$\Ev (\tr( (A(s,t)+\delta Q)^{-1}\partial_t A)$ is smooth on $U$,
and 
$$\partial_s \Ev \big(\tr( (A(s,t)+\delta Q)^{-1}\partial_t A)\big)
= \partial_t \Ev \big(\tr( (A(s,t)+\delta Q)^{-1}\partial_s A)\big).$$
\end{prop}

This proposition has the following corollary.
\begin{corollary}\label{cor:1}
Let $M$, $B$, be $d\times d$ matrices independent of $s$ and $t$.  Then
\begin{align*} & 
\int_0^1 \Ev \Big( \tr  \big( 
(s {B}+{M}+\delta {Q})^{-1} {B}\big) \Big)ds 
= \int_0^1 \Ev
\Big( \tr  \big( 
( {B}+t {M}+\delta {Q})^{-1} {M}\big) \Big)dt
\\ & \ \ \ - \, \int_0^1 \Ev \Big( \tr  \big( 
( t {M}+\delta {Q})^{-1} {M}\big) \Big)dt
\hspace{4mm} +  \int_0^1 \Ev \Big( \tr  \big( 
( s{B}+\delta {Q})^{-1} {B}\big) \Big)ds.
\end{align*}
\end{corollary}
\begin{proof}
Using the previous proposition, this follows from the Fundamental
Theorem of Calculus:
\begin{align*} & 
\int_0^1 \Ev \Big( \tr  \big( 
(s {B}+{M}+\delta {Q})^{-1} {B}\big) \Big)ds -\int_0^1 \Ev \Big( \tr
  \big( 
( s{B}+\delta {Q})^{-1} {B}\big) \Big)ds
\\ & 
= \int_0^1 \partial_t \int_0^1  \Ev \tr\left( 
(s {B}+t{M}+\delta {Q})^{-1}B\right)  ds dt\\ & 
= \int_0^1 \partial_s \int_0^1 \Ev \tr \left(  
(s {B}+t{M}+\delta {Q})^{-1}M\right) dt ds
\\ & = \int_0^1 \Ev
\Big( \tr  \big( 
( {B}+t {M}+\delta {Q})^{-1} {M}\big) \Big)dt
-  \int_0^1 \Ev \Big( \tr  \big( 
( t {M}+\delta {Q})^{-1} {M}\big) \Big)dt.
\end{align*}
\end{proof}

Proposition \ref{p:trick} follows from the subsequent two lemmas.
\begin{lem}
Let $A(s,t)$, $B(s,t)$ be $d\times d$ matrices depending smoothly on 
$(s,t)\in U\subset \Complex^2$.  With $Q$ a random matrix as in Proposition
\ref{p:trick} and $\delta>0$, 
$$\Ev\left( \tr\big( (A(s,t)+\delta Q)^{-1} B(s,t)\big)\right)
\in C^{\infty}(U).$$
\end{lem}
\begin{proof}
We prove the lemma by writing the expected value
as an integral:
\begin{align*}
\Ev \left( \tr( (A+\delta Q)^{-1} B)\right)
& = \int 
\tr( (A+\delta Q)^{-1} B) e^{-\| Q\|^2_{HS}} d\mcl(Q)\\
& = \int \tr( (\delta Q)^{-1} B) e^{-\| Q-\frac{1}{\delta}A \|^2_{HS}} 
d\mcl(Q).
\end{align*}
Now, for a $d\times d $ matrix $\tilde{B}$, 
$|\tr( (\delta Q)^{-1} \tilde{B})| \leq {C}
 | \det Q|^{-1} \|\tilde{B}\|
 \|Q\|^{d-1} / {\delta} $, 
where the constant $C$ depends on $d$.
Moreover,
$$|\partial_s^j \partial_t^k e^{-\|Q-\frac{1}{\delta} A\|^2_{HS}}
| \leq C_{j,k,d}(\sum_{j'\leq j, k'\leq k}\| \partial_s^{j'}\partial_t^{k'}
A\|)\left( \frac{\|Q\|}{\delta^2}\right)^{j+k}e^{-\|Q-\frac{1}{\delta} A\|^2_{HS}}.$$
Since, using Lemma \ref{l:L1} 
$\int  | \det Q|^{-1} 
 (1+\|Q\|)^{m} e^{-\| Q-\frac{1}{\delta} A \|^2_{HS}}d\mcl (Q)<\infty$,
for any finite $m$, the smoothness of $A$ and $B$ proves the lemma.
\end{proof}

If $M$ is an invertible matrix depending smoothly on $s$ and $t$,
then 
\begin{equation}\label{eq:logder}
\tr( M^{-1} \partial_t M)= \frac{\partial _t \det M}{\det M} \; \text{and}\; 
\partial_ s \tr( M^{-1} M_t)= \partial _t \tr( M^{-1}M_s) .
\end{equation}
The lemma below shows that something similar is true when taking expected
values, even though the matrices under consideration are not  
invertible for some values of the random variable.
\begin{lem}\label{eq:switchorder}
Let $A(s,t)$ be a $d\times d$ matrix depending smoothly on 
$(s,t)\in U\subset \Complex^2$, and $Q$ a random matrix as in Proposition
\ref{p:trick}.  Then for $\delta >0$
$$\partial_s \Ev \left( \tr( (A+\delta Q)^{-1} \partial_t A)\right)
= \partial_t \Ev \left( \tr( (A+\delta Q)^{-1} \partial_s A)\right).$$
\end{lem}
\begin{proof}
Let $\chi_{\epsilon}\in C^{\infty}(\Real)$ satisfy 
$\chi_{\epsilon}(x)=1$ for $|x|<\epsilon/2$ and $\chi_{\epsilon}(x)=0$
 for $|x|> \epsilon $.  Then
\begin{align}
\label{eq:chiep}
\nonumber 
\partial_s \Ev \left( \tr( (A+\delta Q)^{-1} \partial_t A)\right)
& = \partial_s \Ev \left( 
\chi_{\epsilon}(\det(A+\delta Q))
 \tr( (A+\delta Q)^{-1} \partial_t A)\right)
\\ & + \partial_s \Ev \left( 
\big(1-\chi_{\epsilon}(\det(A+\delta Q))\big)
 \tr\big( (A+\delta Q)^{-1} \partial_t A\big)\right).
\end{align}
Now 
\begin{align*}& 
\partial_s \Ev \left( 
\big(1-\chi_{\epsilon}(\det(A+\delta Q))\big)
 \tr\big( (A+\delta Q)^{-1} \partial_t A\big)\right)\\ & 
= \int \big(1-\chi_{\epsilon}(\det(A+\delta Q))\big)
\partial_s \tr\big( (A+\delta Q)^{-1} \partial_t A\big)
 e^{-\|Q\|^2_{HS}} d\mcl(Q)\\ & -
\int \chi_{\epsilon}'(\det(A+\delta Q))\big( \partial_s \det (A+\delta Q)\big)
           \tr\big( (A+\delta Q)^{-1} \partial_t A\big)
 e^{-\|Q\|^2_{HS}} d\mcl(Q)
\end{align*}
where we can freely interchange differentiation and 
integration since the integrand is smooth and it and its derivatives
are integrable.
But using (\ref{eq:logder}), we get
\begin{align*}& 
\partial_s \Ev \left( 
\big(1-\chi_{\epsilon}(\det(A+\delta Q))\big)
 \tr\big( (A+\delta Q)^{-1} \partial_t A\big)\right)\\ & 
= \int \big(1-\chi_{\epsilon}(\det(A+\delta Q))\big)
\partial_t \tr\big( (A+\delta Q)^{-1} \partial_s A\big)
 e^{-\|Q\|^2_{HS}} d\mcl(Q)\\ & -
\int \chi_{\epsilon}'(\det(A+\delta Q))\partial_t \det (A+\delta Q)
           \tr\big( (A+\delta Q)^{-1} \partial_s A\big)
 e^{-\|Q\|^2_{HS}} d\mcl(Q)\\ & 
= \partial_t \Ev \left( 
\big(1-\chi_{\epsilon}(\det(A+\delta Q))\big)
 \tr\big( (A+\delta Q)^{-1} \partial_s A\big)\right).
\end{align*}

On the other hand, the first term on the right in (\ref{eq:chiep})
satisfies
\begin{align*} & 
\lim_{\epsilon \downarrow 0}
\partial_s \Ev \left( 
\chi_{\epsilon}(\det(A+\delta Q))
 (\tr( (A+\delta Q)^{-1} \partial_t A)\right)\\ &
= 
\lim_{\epsilon \downarrow 0}\partial_s 
\int \chi_{\epsilon}(\det(A+\delta Q))
 (\tr( (A+\delta Q)^{-1} \partial_t A)e^{-\|Q\|^2_{HS}} d\mcl(Q)
\\ & = \lim_{\epsilon \downarrow 0}\partial_s 
\int \chi_{\epsilon}(\det(\delta Q))
 (\tr ( (\delta Q)^{-1} \partial_t A)e^{-\|Q-\frac{1}{\delta}A\|^2_{HS}} 
d\mcl(Q)
= 0
\end{align*}
since $(\tr( (\delta Q)^{-1} \partial_t A)e^{-\|Q-\frac{1}{\delta}A\|^2_{HS}}$
and its $s$ derivative are both in $L^1$, using Lemma \ref{l:L1}.
\end{proof}

\section{Reduction to a deterministic problem}
\label{rdp}

In this section we will show how to reduce the random 
problem problem to a deterministic one. That will be done
using the singular value decomposition of the matrix $ f_N$.

Let $A$ be a square matrix, and let $USV^*$ be a singular
value decomposition for $A$. 
We make the following simple observation: for 
$ \psi \in \CIc ( \RR , \RR ) $ equal to $ 1 $ on $ [-1,1] $,
\begin{equation}
\label{eq:Aal}    ( A + \alpha \psi ( A A ^* / \alpha^2) U V^* )^{-1}
= { {\mathcal O}} ( 1/ \alpha ) \; : \; \ell^2 \longrightarrow \ell^2 \,,
\end{equation} 
which becomes totally transparent by writing
 $ \psi ( A A ^* / \alpha^2) U V ^* = U \psi ( (S / \alpha)^2 ) V^* $.

The random problem is reduced to a deterministic one by 
using an operator of the form \eqref{eq:Aal}.

\begin{prop}
\label{p:ra2de}
For a smooth curve $ \gamma $ define
\begin{equation}
\label{eq:IN}
I_N ( \gamma ) \defeq  \int_\gamma \EE 
\tr ( f_N + \delta Q_N  - z )^{-1}  dz \,
\end{equation}
where $Q_N$ is a complex $N^n \times N^n$ matrix, with entries
indepent $N(0,1)$ random variables.
Let $ f_N = U_N S_N V_N^* $ 
be a singular value decomposition of $ f_N$,
and let  $ \psi \in \CIc ( \RR ; [0,1]) $ be equal to $ 1 $ on $ [-1,1] $.
If  
\begin{equation}
\label{eq:condgaal}
0 \in \gamma \,, \ \ | \gamma | < \alpha/4 \,, \ \ 
\ \  \delta \ll \alpha \,, \end{equation}
then 
\begin{equation}
\label{eq:ra2de}
\begin{split} 
& \int_\gamma \EE \tr ( f_N + \delta Q_N ( \omega ) - z )^{-1} dz  = \, \\
& \ \ \ \ \int_\gamma \EE \tr ({ f_N + \alpha \psi ( 
f_N f_N^*/ \alpha^2 ) U_N V_N^* } + \delta Q_N ( \omega)   - z )^{-1}  dz
+ E_1  = \, \\
& \ \ \ \  \int_\gamma   \tr ( { f_N + \alpha \psi ( 
f_N f_N^*/ \alpha^2 ) U_N V_N^* }  - z )^{-1} dz  
+ E_2 \,, 
\end{split}
\end{equation}
where
\begin{equation}
\label{eq:boundE} E_1, E_2 = {\mathcal O} \left(
 d \log \left( \frac{\alpha} {\delta} \right) + \frac{N^{4n}}{\delta}e^{-\alpha ^2/4(3N^n \delta)^2}
 \right) \,,
\end{equation}
and $ d = \rank \bbbone_{\supp \psi}  \left( { f_N f_N^*} /{\alpha^2} \right) $.
\end{prop}

The proof of this proposition will use the following lemma.
\begin{lem}\label{l:chisides}
Let $f_N,\; U_N, \; S_N,\; V_N,\; \psi,\; \delta,\; d,$  and $\alpha$ be as
in the statement of Proposition \ref{p:ra2de}.  
Let $\chi\in L^{\infty}(\Real)$ be the characteristic function for the 
support of
$\psi $.  Then, if $|z|\leq \alpha/4$,
\[
\left| \int_0^1 \Ev \tr \left( ( f_N +s\alpha \chi ( 
f_N f_N^*/ \alpha^2 ) U_N V_N^*  - z + \delta Q )^{-1} \alpha \chi ( 
f_N f_N^*/ \alpha^2  ) U_N V_N^* \right)   ds
\right|
\]
satisfies the bound \eqref{eq:boundE}.
\end{lem}
\begin{proof}
First suppose that for a $m\times m$ matrix $\tilde A$,
\begin{equation}
\tilde{A}= \left( \begin{array}{cc} \tilde{A}_{11} & \tilde{A}_{12}\\
\tilde{A}_{21} & \tilde{A}_{22} \end{array}\right) \ \ \ \ \
 \text{and} \ \ \  \ \ 
\tilde{A}^{-1}= \left( \begin{array}{cc} \tilde{B}_{11} & \tilde{B}_{12}\\
\tilde{B}_{21} & \tilde{B}_{22} \end{array}\right)
\end{equation}
with $\tilde{A}_{11},$ 
$\tilde{B}_{11}$ $d\times d $ matrices and $\tilde{A}_{11},$ $\tilde{B}_{11}$
$(m-d)\times(m-d)$  matrices.  Then if $\tilde{A}_{22}$ is invertible, 
we have the Schur complement formula, 
\begin{equation}\label{eq:B11}
\tilde{B}_{11}= \left( \tilde{A}_{11} + \tilde{A}_{12}\tilde{A}_{22}^{-1}
\tilde{A}_{21}\right)^{-1} \,,
\end{equation}
see \cite{SjZw07} for a review of some of its applications in spectral theory.

We note, using  $ \psi ( A A ^* / \alpha^2) U V ^* = U \psi ( (S / \alpha)^2 ) V^* $
and the unitarity of $U_N$,
$V_N$, 
\begin{multline}\label{eq:rewrite}  
\Ev \tr \big( (f_N +s\alpha \chi  ( 
f_N f_N^*/ \alpha^2 ) U_N V_N^*  - z+ \delta Q_N )^{-1} \alpha \chi ( 
f_N f_N^*/ \alpha^2 ) U_N V_N^*\big)
\\ 
= \Ev \tr \big( (S_N +s\alpha \chi ( 
S_N S_N^*/ \alpha^2 )   - U_N^* z V_N+ \delta Q_N )^{-1} \alpha \chi ( 
S_N S_N^*/ \alpha^2 ) \big) .
\end{multline}
The main idea of the proof will be to effectively reduce the 
dimension of the matrices we work with, from $N^n$ to $d$.
We can assume that $U_N$, $V_N$, $S_N$ are chosen 
so that the diagonal elements $\sigma_1,...,\sigma_{N^n}$ of $S_N$
satisfy $\sigma_1\leq \sigma_2\cdot\cdot \cdot \leq \sigma_{N^n}$.  
Let $\pr$ denote projection onto the range of $\chi(S_N^2/\alpha^2)$,
which is the same as projection off of the kernel of $\chi(S_N^2/\alpha^2)$.
Then
$$\pr = \left( \begin{array}{cc} I_d & 0 \\ 0 & 0 
\end{array}
\right) .
$$
and $\alpha \chi ( 
S_N^2 / \alpha^2 ) $ takes the form
$$ 
\left( \begin{array}{cc} \alpha I_d & 0 \\ 0 & 0 
\end{array} \right). $$
We also write 
$$ S_N+s \alpha \chi ( 
S_N S_N^*/ \alpha^2 )   - U_N^* z V_N 
= \left( 
\begin{array}{cc} 
s\alpha I_d+A_{11} & A_{12}\\
A_{21} & A_{22}
\end{array}\right), 
$$
 and
$$Q_N = \left( 
\begin{array}{cc} Q_{11  } & Q_{12} \\ Q_{21} & Q_{22} 
\end{array} \right)
$$ where $A_{11},\; Q_{11}$ are  $d\times d$-dimensional matrices, and 
$A_{22},\; Q_{22}$ are $(N^n-d)\times (N^n-d)$-dimensional.  
Since $S_N$ is diagonal and $|z|\leq \alpha/4$, we have 
$\|A_{12}\|\leq \alpha/4$, $\|  A_{21}\|\leq \alpha/4$.

Using this notation, we have that 
$A_{22}$ is invertible, with norm at most  $4/3\alpha$.  
Now restrict $Q_N$ to the set with 
\begin{equation}\label{eq:Qrestrict}
\delta \| Q_N- \pr Q_N \pr\|_{\sup} \leq {\alpha}N^{-n}/4.
\end{equation}
Note that poses no restriction on $Q_{11}$.  For such $Q_N$,
$A_{22}+\delta Q_{22}$ is invertible, with norm at most $2/\alpha$.
Restricting to this set of $Q_N$ and using  
 (\ref{eq:B11}), we find
\begin{gather*}
\tr\Big( (S_N +s\alpha \chi ( 
S_N S_N^*/ \alpha^2 )   - U_N^* z V_N )^{-1} \alpha \chi ( 
S_N S_N^*/ \alpha^2 )\Big) 
= 
\tr_d \left( \alpha (s \alpha I_d +  M_d + \delta Q_{11} )^{-1}\right)\,,
\end{gather*}
where we use the notation $\tr_d$ to emphasize we are taking the 
trace of a $d\times d$ matrix, and  where
$$M_d =  
A_{11}
 - (A_{12}+\delta Q_{12})(A_{22}+\delta Q_{22})^{-1}
(A_{21}+\delta Q_{21})
$$ 
is a $d\times d$ matrix depending on $Q_{12}$,  $Q_{21}$, and $Q_{22}$, 
but not on $Q_{11}$. 
Since $\| A_{11}\| = \| \pr (S_N - zU_N^* V_N)\pr\| \leq C \alpha$
and $\| A\|_{12}\leq \alpha/4$, $\| A \|_{21}\leq \alpha/4$, we have 
$ \| M_d\| \leq C \alpha$, for a new constant $C$ independent of $N$, $\delta$,
 and 
$Q_N$ satisfying (\ref{eq:Qrestrict}).
 
Next we take the expected value in the $Q_{11}$ variables only:
$$\Ev_{Q_{11}}(F(Q_N))= \int_{Q_{11}\in \Complex^{d^2 }} F(Q_N)
 e^{-\|Q_{11}\|_{HS}^2}
d{\mathcal L}(Q_{11}).$$  Still requiring $Q_N$
 to satisfy (\ref{eq:Qrestrict}),
which is not a restriction on $Q_{11}$, 
and using Corollary \ref{cor:1}, 
we get
\begin{align*}& 
\int_0^1 \Ev_{Q_{11}}\Big(\alpha 
 \tr_d \big( M_d+s \alpha I_d+\delta Q_{11})^{-1}\big) 
\Big)ds = \,  \\
& \ \ \ \int_0^1 \Ev_{Q_{11}}
\Big( \tr_d \big( tM_d+ \alpha I_d +\delta Q_{11})^{-1}M_d\big)
\Big)dt
- \int_0^1 \Ev_{Q_{11}}\Big( \tr_d \big( s \alpha I_d +\delta Q_{11})^{-1}
\alpha \big)\Big) ds
\\ & \ \ \ \ \ \
+ \int_0^1 \Ev_{Q_{11}}\Big( \tr_d \big( t M_d+ \delta Q_{11})^{-1} M_d 
\big)\Big) dt.
\end{align*}
 
Recalling that 
$ \|M_d\| \leq C \alpha$ we see from Proposition
\ref{p:tA} that the second and third terms on the right are 
${\mathcal O}(d \log (\alpha/\delta))$, if $\alpha/\delta >e$.
Moreover, 
$$\| M_d - \pr S_N \pr\| \leq \frac{\alpha}{2},$$
and $S_N \geq 0$.  Therefore,
                   for $0\leq t\leq 1$, $\alpha I_d +t M_d$ 
is invertible, with the 
inverse having norm at most $3/\alpha $.  Thus from Lemma 
\ref{l:invertible} we see that 
$$ \left |\int_0^1 \Ev_{Q_{11}}
\Big( \tr_d \big( t M_d + \alpha I_d +\delta Q_{11})^{-1} M_d \big)
\Big)dt \right| = {\mathcal O}( d) 
+ {\mathcal O}\left( \frac{d^4}{\delta} e^{-\alpha^2/4 (3 d \delta)^2}\right).$$
The implicit constants in both cases are independent of $Q-\pr Q\pr$ 
satisfying (\ref{eq:Qrestrict}).  Thus we get
\begin{gather}\label{eq:firstpart}
\begin{gathered}
\int_0^1 
\Ev\Big( \tr \big((S_N+s\alpha \chi(S_N^2/\alpha^2)  +\delta Q-z U_N ^* V_N)
^{-1} 
\alpha \chi(S_N^2/\alpha^2) \big)
\bbbone_{\{ \delta \| Q-\pr Q \pr\|_{\sup}\leq \frac{\alpha}{4 N^n }
\}}
 \Big) ds \\ =   {\mathcal O}( d\log (\alpha/\delta)) 
+ {\mathcal O}\left( {d^4} {\delta}^{-1} e^{-\alpha^2 /4 (d 3 \delta)^2}\right)
\end{gathered}
\end{gather}
where for a set $E$, $\bbbone_E$ is the characteristic function of $E$.

Exactly as in the proof of the Lemma \ref{l:invertible}, we can 
show that 
\begin{gather}\label{eq:part2}
\begin{gathered}
\Ev\Big( \tr \big((S_N+s\alpha \chi(S^2_N\alpha^2) +\delta Q_N-z U_N^*V_N )^{-1} 
\alpha \chi(S^2_N\alpha^2)\big)
\bbbone_{\{ \delta \| Q-\pr Q \pr\|_{\sup}> \alpha/(4N^n)\}} \Big)
 \\
= {\mathcal O} \left( {N^{4n}} {\delta}^{-1} e^{-\alpha^2/(4 N^n\delta)^2}\right).
\end{gathered}
\end{gather}
Using (\ref{eq:rewrite}), 
(\ref{eq:firstpart}), and (\ref{eq:part2}), we prove the lemma.
\end{proof}

We now use Lemma \ref{l:chisides} in a preliminary step towards
 proving Proposition \ref{p:ra2de}.
\begin{lem}\label{l:withchi}
Let $f_N,\; U_N, \; S_N,\; V_N,\; \psi,\; \delta,\; d,$  and $\alpha$ be as
in the statement of Proposition \ref{p:ra2de}, and set $\chi=\bbbone_{\supp \psi}$.
Then 
\[ \begin{split}
& \int_{\gamma} \EE
\tr ( f_N + \delta Q_N  - z )^{-1}  dz 
= \int_\gamma \EE \tr ({ f_N + \alpha \chi ( 
f_N f_N^*/ \alpha^2 ) U_N V_N^* } + \delta Q_N   - z )^{-1}  dz \\ 
& \ \ \  + {\mathcal O}\left(d \log \left(\frac{\alpha}{\delta}\right)\right)
+ {\mathcal O}\left( \frac{N^{4n}}{\delta}e^{-\alpha ^2/4(3N^n \delta)^2}\right).
\end{split}
\]
\end{lem}
\begin{proof}  The proof uses the same type of argument as Corollary
\ref{cor:1}.  
Using the Fundamental Theorem of Calculus, 
\begin{align*} & 
 \int_{\gamma} \EE
\tr ( f_N + \delta Q_N  - z )^{-1}  dz
- \int_\gamma \EE \tr ({ f_N + \alpha \chi ( 
f_N f_N^*/ \alpha^2 ) U_N V_N^* }+ \delta Q_N   - z )^{-1}  dz \\ & 
= -\int_0^1 \partial_s \int_\gamma \EE  \tr ({ f_N + s \alpha \chi ( 
f_N f_N^*/ \alpha^2 ) U_N V_N^* }+ \delta Q_N   - z )^{-1}
   dz ds\\ &
= -\int_{\gamma} \partial_z \int_0^1 \EE \left( 
 \tr ({ f_N + s \alpha \chi ( 
f_N f_N^*/ \alpha^2 ) U_N V_N^* }+ \delta Q_N   - z )^{-1}  \alpha \chi(f_Nf_N^*/\alpha^2)\right) ds dz 
\end{align*}
where we use Proposition \ref{p:trick}.  The right hand side is
$$\sum_{\pm} \mp \int_0^1 \EE \tr \left( ({ f_N + s \alpha \chi ( 
f_N f_N^*/ \alpha^2 ) U_N V_N^* }+ \delta Q_N   - z_{\pm} )^{-1}
\alpha \chi(f_Nf_N^*/\alpha^2)\right)   ds
$$
where $z_{\pm}$  are the endpoints of $\gamma$.  Then using Lemma \ref{l:chisides}
finishes the proof.
\end{proof}

We are now able to give a straightforward proof of Proposition \ref{p:ra2de}.
\begin{proof}[Proof of Proposition \ref{p:ra2de}]
We begin by noting that, with $\chi= \bbbone_{\supp \psi}$
$$\| ({ f_N + \alpha \chi ( 
f_N f_N^*/ \alpha^2 ) U_N V_N^* }   - z )^{-1} \| = {\mathcal O}(1/\alpha)$$
and 
$$\|({ f_N + \alpha \psi( 
f_N f_N^*/ \alpha^2 ) U_N V_N^* }   - z )^{-1} \| = {\mathcal O}(1/\alpha)$$
when $|z| \leq \alpha/4$.
Moreover, the rank of $\chi(f_N f_N^*/ \alpha^2 )$ is $d$ and the rank of 
 $\psi (f_N f_N^*/ \alpha^2 )$ is at most $d$, and both operators
have norm at most $1$.  Then
\begin{align*} & \left| \int_\gamma 
\left( \tr ({ f_N + \alpha \chi ( 
f_N f_N^*/ \alpha^2 ) U_N V_N^* }  - z )^{-1} 
- \tr ({ f_N + \alpha \psi ( 
f_N f_N^*/ \alpha^2 ) U_N V_N^* }  - z )^{-1}\right) 
 dz \right|\\
& = \alpha \left| \int_\gamma \tr\left( ({ f_N + \alpha \chi ( 
f_N f_N^*/ \alpha^2 ) U_N V_N^* }  - z )^{-1} \left( 
\chi(f_N f_N^*/ \alpha^2 ) - \psi ( 
f_N f_N^*/ \alpha^2 )\right)\right.\right. \\ & \hspace{40mm}\left. \left. \times 
({ f_N + \alpha \psi ( 
f_N f_N^*/ \alpha^2 ) U_N V_N^* }  - z )^{-1} \right) dz \right|
\\ & \leq \int_{\gamma} \frac{C d}{\alpha} dz = {\mathcal O}(d).
\end{align*}
Thus, applying 
Lemmas \ref{l:withchi} and  \ref{l:invertible} proves the Proposition.
\end{proof}

\section{Proof of Theorem}
\label{pot}

The proof of Theorem will be deduced from the following 
local result:

\begin{prop}
\label{p:loc}
Under the assumption of the main theorem, 
let $ \gamma \subset \partial \Omega $ be a connected segment of 
length 
\begin{equation}
\label{eq:lgamma}  
| \gamma | \leq  \frac  \alpha  C \,, \ \
 h = \frac 1 { 2 \pi N } \,, \ \ \alpha = h^{\rho} \,, \ \ 0 < \rho < \frac 12 
\,
\end{equation}
and let $I_N(\gamma)$ be as defined by (\ref{eq:IN}).
Then for
$\exp \left( { - h^{-\epsilon} } \right) < \delta < h^{p_0} $,
we have 
\begin{equation}
\label{eq:IN1}
I_N ( \gamma ) = N^n \int_\gamma \int_{ \TT^{2n} } ( f ( w) - z )^{-1} 
d {\mathcal L} ( w) d z + {\mathcal O} ( |\gamma| 
h^{-n+ \rho( \kappa - 1 ) - 2 \epsilon }  )
+ {\mathcal O} ( |\gamma |  h^{-n + 1 - 2 \rho } ) 
\,,
\end{equation}
where we note that \eqref{eq:cond} with $ \kappa > 1 $ implies
that $ ( f ( w ) - z ) ^{-1} \in L^1 ( \TT^{2n} ) $ so that the first
term on the right hand side makes sense.
\end{prop}

Assuming the proposition we easily give the 

\noindent
{\em Proof of Theorem.} 
We divide $ \partial \Omega $ into 
$ J= C'/ \alpha $ disjoint 
segments $ \gamma_j $, 
$ | \gamma_j |  \leq \alpha/ C$. Proposition 
\ref{p:loc} implies that 
\begin{gather*}
\EE \left( \tr \int_{\partial  \Omega } ( f_N + \delta Q_N - z ) ^{-1} \right) 
 = \sum_{j=1}^J I_N ( \gamma_j ) = \, \\
N^n \int_{\partial \Omega } \int_{ \TT^{2n} } ( f ( w) - z )^{-1} 
d {\mathcal L} ( w) d z 
 +  {\mathcal O} (
h^{-n+ \rho({\kappa - 1 }) -2 \epsilon  })  + {\mathcal O} (  h^{-n + 1 - 2 \rho } ) 
\,. \end{gather*}
We now choose $ \rho = 1/( \kappa + 1 ) $, to optimize the 
error, that is to arrange,
$ \rho ( \kappa - 1 ) = 1 - 2 \rho $. That means that the error is
$ {\mathcal O} (N^{n-\beta } ) $ for any 
$ \beta < 1 - 2 \rho = ( \kappa - 1 )/({ \kappa + 1 }) $.

Hence 
\[ \begin{split}  {\mathbb E}_\omega \left( | {\rm{Spec}}\, 
( f_N + N^{-p} Q_N( \omega) ) \cap \Omega | \right)  & =  
 \frac{ 1 } { 2 \pi i } \int_{\partial \Omega} 
\EE \tr ( f_N + N^{-p } Q_N ( \omega ) - z )^{-1} dz  \\
 & = 
\frac{ 1 } { 2 \pi i } 
\int_{\partial \Omega}   N^n \int_{\TT^{2n} } \frac{ d {\mathcal L} ( w )} 
{  { f ( w ) - z }     } dz + {\mathcal O} ( N^{n-\beta} ) \\
& = 
   N^n \int_{\TT^{2n} } \frac{ 1 } { 2 \pi i } 
\int_{\partial \Omega} \frac{ dz }   {  { f ( w ) - z }     }
 d {\mathcal L} ( w ) + {\mathcal O} ( N^{n-\beta} )  \\
&  =  {  N^n {\rm{vol}}_{ \TT^{2n} } ( f^{-1} ( \Omega ) ) } + 
{\mathcal O} ( N^{n-\beta} ) \,, 
\end{split}\]
which is the statement of the theorem.
\stopthm

\medskip

\noindent
{\em Proof of Proposition \ref{p:loc}.} 
Without loss of generality we can assume that $ 0 \in \gamma $.
From Proposition \ref{p:ra2de}
we already know that $ I_N ( \gamma ) $ can be approximated by 
a deterministic expression 
\begin{equation}
\label{eq:defIt} \widetilde I_N ( \gamma) \defeq 
\int_\gamma   \tr ( { f_N + \alpha \psi ( 
f_N f_N^*/ \alpha^2 ) U_N V_N^* }  - z )^{-1} dz  \,,
\end{equation}
with, if $\alpha/\delta N^n \gg 0$
\[ I_N ( \gamma ) - \tilde I_N (\gamma ) = 
 {\mathcal O} \left( 
 e^{-c_o \alpha/N^n\delta} +  d \log \left( \frac{\alpha} {\delta} \right) 
\right) \,, \]
for some $c_0>0$,
where $ d $ is the rank of $ \psi (f_N f_N^*/ \alpha^2 ) $. 
We choose $ \alpha $ as in \eqref{eq:lgamma}, $ \alpha = h^{\rho }$, where 
\[ h = \frac 1 { 2 \pi N} \,, \ \ 0 <  \rho < \frac 1 2 \,. \]
In view of Proposition \ref{p:rank}, $ d = {\mathcal O} ( 
h^{-n + \rho \kappa } )  $ and this shows that for this choice of 
$ \alpha $ and for $ \delta $ satisfying the condition in the proposition,
with $ p_0 > n+1/2$,
\[ I_N ( \gamma ) - \tilde I_N (\gamma ) = 
{\mathcal O} ( h^{-n+  \epsilon + \kappa \rho } + \exp(-c_0h^{ n - p_0 + \rho } ) =
 {\mathcal O} (|\gamma|  h^{ -n + ( \kappa -  1 ) \rho -  \epsilon } ) \,. \]


Thus we will prove \eqref{eq:IN1} by showing that
\begin{equation}
\label{eq:cru}
\begin{split} 
\tr ( f_N + \alpha \psi ( 
f_N f_N^*/ \alpha^2 ) U_N V_N^*  - z )^{-1}  & = N^n 
\int_{ \TT^{2n} } ( f ( w) - z )^{-1} 
d {\mathcal L} ( w)  \\
& \ \ \ \ \ \ \ \ \ + {\mathcal O} (  h^{-n + 1 - 2 \rho } ) + 
{\mathcal O}( h^{-n+ \rho({\kappa - 1} ) } ) \,. 
\end{split}
\end{equation}
We first show that it is enough to consider $ z = 0 $. In fact,
let $ U_N ( z) S_N ( z ) V_N ( z ) ^* $ be the singular
value decomposition of $ f _N - z $, and put
\[  B_N ( z , w ) \defeq 
 ( f_N - w + \alpha \psi ( (f_N - z ) ( f_N - z ) ^*/ \alpha^2 ) U_N ( z ) 
V_N^* ( z ) )^{-1} \,.\]
Then $ \tr\left( B_N ( z , z  ) - B_N ( 0 , z ) \right)= \, $ 
\begin{gather*} 
\alpha \tr\left( B_N ( 0, z ) 
\left( \psi ( 
f_N f_N^*/ \alpha^2 ) U_N V_N^* - \psi ((f_N - z ) ( f_N - z ) ^*/ \alpha^2 ) ) U_N ( z ) V_N^*( z)   \right) 
B_N ( z, z ) \right).
\end{gather*}
Since $ \rank \psi ((f_N - z ) ( f_N - z ) ^*/ \alpha^2 ) ) = 
{ \mathcal O }  ( h^{-n+\kappa \rho  } ) $ for $ z \in \gamma $, 
and $ B ( z , w ) = {\mathcal O}_{\ell^2 \rightarrow 
\ell^2 }  ( 1/\alpha ) $ for $ |z - w| \leq \alpha / C' $, 
we obtain
\[  \tr\left( B_N ( z , z  ) - B_N ( 0 , z ) \right) = 
{\mathcal O}  (  h^{-n + \rho ( {\kappa-1} )  }  )  \,,\]
which can be absorbed in the error on the right hand side of \eqref{eq:cru}.
Thus we only need to prove \eqref{eq:cru} with the left hand side
replaced by $ B( z, z ) $ and we can simply take $ z = 0 $.

In other words we now want to prove 
\begin{equation}
\label{eq:cru0}
\begin{split} 
\tr ( f_N + \alpha \psi ( 
f_N f_N^*/ \alpha^2 ) U_N V_N^*  )^{-1}  & = N^n 
\int_{ \TT^{2n} } \frac{d {\mathcal L} ( w) } { f ( w) }
+ {\mathcal O} (  h^{-n + 1 - 2 \rho } ) + 
{\mathcal O}( h^{-n+ \rho({\kappa - 1} ) } ) \,. 
\end{split}
\end{equation}

The difficulty lies in the fact that the operators $
 f_N + \alpha \psi ( 
f_N f_N^*/ \alpha^2 ) U_N V_N^*   $ do {\em not} seem to have a nice
microlocal characterization. We are 
helped by the following identity: if $ \tilde \psi \in \CIc ( \RR , [0,1] ) $ 
is equal to $ 1 $ on the support of
$ \psi $ then 
\begin{gather}
\label{eq:svd1}
\begin{gathered}
  ( 1 - \tilde \psi ( f_N^* f_N / \alpha^2) )   ( f_N + \alpha \psi ( 
f_N f_N^*/ \alpha^2 ) U_N V_N^* )^{-1} =  \\
\ \ \ \ \ \ \ ( 1 - \tilde \psi ( f_N^* f_N / \alpha^2) )  
f_N^* ( f_N f_N^* + \alpha^2 \psi ( f_N f_N^*/ \alpha^2 )  )^{-1} \,.
\end{gathered}
\end{gather}
This is a consequence of an identity from linear algebra:

\begin{lem}
\label{l:la}
Let $ A $ be a matrix and $ U S V^* $ be its singular value
decomposition. If $ \psi, \tilde \psi \in \CIc ( \RR ; [0,1] ) $, 
$ \psi $ is equal to $ 1 $ on $ [ -1, 1] $, and $ \tilde \psi $ is 
equal to $ 1 $ on the support of $ \psi $, then 
\begin{equation}
\label{eq:la}
  ( 1 - \tilde \psi ( A^* A ) ) ( A + 
\psi ( A A^* ) U V^* )^{-1}  =    ( 1 - \tilde \psi ( A^* A ) )  A^* ( A  A^* + 
\psi ( A A^* )  )^{-1}  \,. 
\end{equation}
\end{lem}
\begin{proof} We first note that 
\[  A^* A = V S^2 V^* \,, \ \ \tilde \psi ( A^* A ) = V \tilde 
\psi ( S^2 ) V^* \,, \]
and similarly $ \psi ( A A^* ) = U \psi ( S^2 ) U^* $. Since
$ S $ is a diagonal matrix, and $ ( 1 - \tilde \psi ) \psi \equiv 0 $,
we get
\[ \begin{split}   ( 1 - \tilde \psi ( A^* A ) ) ( A + 
\psi ( A A^* ) U V^* )^{-1} & = 
V ( 1 - \tilde \psi ( S^2 ) ) V^* V ( S + \psi ( S^2 ) )^{-1} U^* \\
& = V ( 1 - \tilde \psi ( S^2 ) ) ( S + \psi ( S^2 ) )^{-1}  U^*  \\
& = V ( 1 - \tilde \psi ( S^2 ) ) S ( S^2 + \psi ( S^2 ))^{-1} U^* \\
& = \left( V ( 1 - \tilde \psi ( S^2 ) ) V^* \right) 
\left( V S U^*\right) \left( U 
( S^2 + \psi ( S^2 ) )^{-1} U^* \right) \\
& =    ( 1 - \tilde \psi ( A^* A ) )  A^* ( A  A^* + 
\psi ( A A^* )  )^{-1}  \,, \end{split}\]
concluding the proof. 
\end{proof}

The identity \eqref{eq:svd1} follows from \eqref{eq:la} by putting 
$ A = f_N/\alpha $, $ U = U_N $, and $ V = V_N $. Using this we 
we will find a new expression for the left hand side of 
\eqref{eq:cru0} so that the identification with the right hand
side will follow from a suitable semiclassical operator calculus.

\begin{lem}
\label{l:svd2psu}
We have the following approximation for the left hand side of 
\eqref{eq:cru0}:
\begin{equation}
\label{eq:svd2psu}
\tr ( f_N + \alpha \psi ( 
f_N f_N^*/ \alpha^2 ) U_N V_N^*  )^{-1}  = 
\tr f_N^* ( f_N f_N^*+ \alpha^2  \psi ( 
f_N f_N^*/ \alpha^2 )   )^{-1} +  {\mathcal O}  ( h^{-n+\rho({\kappa - 1} ) } ) \,. 
\end{equation}
\end{lem}
\begin{proof}
We use \eqref{eq:svd1} and 
first note that $ 1- \tilde \psi $ can be removed from the left hand
side since
\begin{equation}
\label{eq:err1}
\begin{split} 
& \tr  \tilde \psi ( f_N^* f_N / \alpha^2)  ( f_N + \alpha \psi ( 
f_N f_N^*/ \alpha^2 ) U_N V_N^* )^{-1} = \\  
& \ \ \ {\mathcal O} \left( (\rank
\tilde \psi \left( { f_N f_N^*}/ {\alpha^2} \right) ) 
\| ( f_N + \alpha \psi ( f_N f_N^*/ \alpha^2 ) U_N V_N^* )^{-1} \| \right) 
= {\mathcal O} \left(  {h^{-n+\rho({ \kappa - 1})} } \right) \,.
\end{split}
\end{equation}
The same argument works for the right hand side once we observe that
\[    f_N^* ( f_N f_N^*+ \alpha^2  \psi ( 
f_N f_N^*/ \alpha^2 )   )^{-1} = {\mathcal O}_{\ell^2 \rightarrow \ell^2 }
( 1 / \alpha ) \,, \]
and this follows from using the singular value decomposition since
for non-negative diagonal matrices 
\[  S_N ( S_N^2 + \alpha^2 \psi ( S_N^2 / \alpha^2 ) )^{-1} \leq 1 / \alpha 
\,.\]
\end{proof}

In view of \eqref{eq:cru0} and the lemma we have to prove 
\begin{equation}
\label{eq:cru1}
\tr f_N^* ( f_N f_N^*+ \alpha^2  \psi ( 
f_N f_N^*/ \alpha^2 )   )^{-1} = 
 N^n 
\int_{ \TT^{2n} } \frac{d {\mathcal L} ( w) } { f ( w) }
+ {\mathcal O} ( h^{-n + 1 -2 \rho } )
 + 
{\mathcal O} ( h^{-n+ ( {\kappa - 1} ) \rho } ) 
\,,
\end{equation}
but that follows from the calculus developed in \S \ref{qot}.
In fact, with the $ \alpha $-order function 
$ m ( w , \alpha ) = \alpha^2 + | f( w )|^2 $, given in 
Lemma \ref{l:ord}, 
\begin{gather*}
  f_N f_N^* + \alpha^2  \psi ( f_N f_N^*/ \alpha^2 )  = T_N \,, \ \
T \in S ( m , \alpha ) \,,\\
T = T_0 + h^{1 - 2 \rho } T_1 \,, \ \ T_0 ( w ) = 
| f( w ) |^2 + \alpha^2 \psi ( | f ( w ) |^2/ \alpha^2 ) \,, \ \ 
T_1 \in S ( m , \alpha ) \,. 
\end{gather*}
where we also applied Lemma \ref{l:func}. We also have $ T_0 \geq m/2$
and hence
\[  1 /{ T_0} \in S ( 1 / m , \alpha ) \,,  \ \  1/T \in S (1/m, \alpha ) \,. \] 
Since $ f \in S ( \sqrt  m  , \alpha ) $, we conclude that
\begin{gather*}   f_N^* ( f_N f_N^*+ \alpha^2  \psi ( 
f_N f_N^*/ \alpha^2 )   )^{-1}  = P_N \,, \ \ 
P \in S ( 1/\sqrt  m , \alpha ) \,, \\
P = P_0 + h^{1 - 2 \rho } P_1 \,, \ \ P_1 \in S ( 1/\sqrt  m ) \,, \ \
P_0 ( w ) = \frac{\bar f ( w )}{ | f ( w ) |^2 + \alpha^2 \psi
( |f ( w )|^2 /\alpha^2 ) } \,.
\end{gather*}
We now apply Lemma \ref{l:trace} and obtain (with $ k \gg n $)
\[ \begin{split} 
& \tr f_N^* ( f_N f_N^*+ \alpha^2  \psi ( 
f_N f_N^*/ \alpha^2 )   )^{-1}  = 
 N^n \int_{ \TT^{2n} } P ( w ) d{\mathcal L} ( w ) + 
{\mathcal O} (N^{-k+n} ) \sup_{ |\beta | \leq k } \int | \partial^\beta P  |
d{\mathcal L}  \\ 
&  \ \ \ = \,  N^n \int_{ \TT^{2n} } P_0 ( w ) d{\mathcal L} ( w ) + 
{\mathcal O} ( h^{-n+( 1 - 2 \rho ) } 
+ h^{-n+k( 1 -\rho) } ) 
\int_{ \tn } m ( w , \alpha ) ^{-1/2 }  
 d{\mathcal L} ( w ) 
\end{split}
\]
We have $ m(w, \alpha)^{-1/2} \leq |f( w ) |^{-1} $ and \eqref{eq:cond}
at $ z = 0 $ with $ \kappa > 1 $ implies that $ |f( w ) |^{-1}  $ is 
integrable ($ \kappa = 1 $ would mean that $ |f( w ) |^{-1} $ is in weak
$ L^1 $):
\[ \int_{\tn }   |f( w ) |^{-1}  d {\mathcal L} ( w ) =
\int_0^\infty {\mathcal L} (\{ |f ( w) | < t \} ) t^{-2} dt 
= \int_0^\infty {\mathcal O} ( \min ( t^\kappa , 1 ) ) t^{-2 } dt
< \infty \,. \]
It remains to show that
\begin{equation}
\label{eq:remains}
 \int_{ \TT^{2n} } |  P_0 ( w ) - f ( w ) ^{-1} | d{\mathcal L} ( w ) 
=  {\mathcal O} ( h^{\rho( \kappa -1 ) })  \,. 
\end{equation}
Putting $ \varphi ( x ) \defeq \psi ( x^2) $, we rewrite the left hand 
side above as 
\[ \begin{split} 
& \int_0^\infty {\mathcal L} (\{ |f ( w) | < t \} )  
\partial_t \left( \frac { - \alpha^2 \varphi ( t/\alpha ) } 
{ t ( t^2 + \alpha^2 \varphi ( t/ \alpha ) } \right) dt  = 
\int_0^\infty  {\mathcal L} (\{ |f ( w) | < t \} )   
 \frac{ - \alpha^2 (t/\alpha)  \varphi' ( t/\alpha ) } { t^2 ( t^2 + \alpha^2 
\varphi ( t/ \alpha ))}    dt \\
& \ \ \ \ \ \ \ \ \ \ \ + \, \int_0^\infty {\mathcal L} (\{ |f ( w) | < t \} )  
\frac{\alpha^2 \varphi( t/\alpha) 
(3 t^2 + \alpha^2 \varphi( t/\alpha ) + \alpha^2 (t /\alpha) 
\varphi'( t/\alpha) )}
{ t^2 ( t^2 + \alpha^2 \varphi ( t/\alpha ) )^2}  dt \\
& \ \ \ \ \ \ \ \ \ \ \ \ \ \ \ \ \ \ \ \ \ \ \ \ \ 
\ \ \ \ \ \ \ \ \ \ \ 
\leq \, C \int_0^{2 \alpha} t^{\kappa -2}  dt = C' \alpha^{\kappa - 1 } \,,
\end{split} \]
which is \eqref{eq:remains}. Since we have now established 
\eqref{eq:cru1} this also completes the proof of Proposition \ref{p:loc}.
\stopthm


\begin{thebibliography}{99}

\bibitem{BM} E. Bierstone and P.D. Milman, 
{\em  Semianalytic and subanalytic sets,} 
Publ.~IH\'ES, {\bf 67}(1988), 5--42.

\bibitem{Bo} W.~Bordeaux Montrieux, personal communication.

\bibitem{BoSj}
W.~Bordeaux Montrieux and J.~Sj\"ostrand,
{\em Almost sure Weyl asymptotics for non-self-adjoint elliptic operators on compact manifolds}, {\tt arXiv:0903.2937}.

\bibitem{Bor} 
D. Borthwick, {\em Introduction to K\"ahler quantization},
 In {\em First Summer school in analysis and
mathematical physics}, 
Cuernavaca, Mexico. Contemporary Mathematics series {\bf 260}, AMS(2000), 91--132.

\bibitem{BoU} 
D. Borthwick and A. Uribe, 
{\em On the pseudospectra of Berezin-Toeplitz operators}, 
Methods and Applications of Analysis {\bf 10} (2003), 31--65.

\bibitem{BouzDB}
A.~Bouzouina  and S.~De~Bi\`{e}vre, {\it Equipartition of the eigenfunctions 
of quantized                                                                  
ergodic maps on the torus}, Commun. Math. Phys. {\bf 178} (1996) 83--105.

\bibitem{ChTr} 
S.J.~Chapman and L.N.~Trefethen,
{\em Wave packet pseudomodes of twisted Toeplitz matrices,}
Comm. Pure Appl. Math. {\bf 57} (2004), 1233--1264.

\bibitem{DeSjZw} 
N. Dencker, J. Sj\"ostrand, and M. Zworski, 
{\em Pseudospectra of semi-classical (pseudo)differential operators,}
Comm. Pure Appl. Math., {\bf 57} (2004), 384--415.

\bibitem{DiSj} M. Dimassi and J. Sj\"ostrand, {\it Spectral Asymptotics in
the semi-classical limit,} Cambridge University Press, 1999.

\bibitem{EZ} L.C. Evans and M. Zworski, {\em Lectures on Semiclassical 
Analysis,}\\
{\tt http://math.berkeley.edu/$\sim$zworski/semiclassical.pdf}

\bibitem{FSV} P.~Forrester, N.~Snaith, and V.~Verbaarschot, {\em Introduction
Review} to {\em Special Issue on Random Matrix Theory},
Jour. Physics A: Mathematical and General {\bf36} (2003), R1–-R10.


\bibitem{Ha} M. Hager, unpublished, 2007.

\bibitem{h-sj} 
M. Hager and J. Sj\"ostrand,
{\em Eigenvalue asymptotics for randomly perturbed non-selfadjoint operators.}
Math. Ann. {\bf 342} (2008), no. 1, 177--243. 

\bibitem{Me} N.~Mehta, {\em Random Matrices}, 3rd Edition, Elsevier, 
Amsterdam, 2004.

\bibitem{NZ} S. Nonnenmacher and M. Zworski, {\em Distribution
of resonances for open quantum maps,}
Comm. Math. Phys. {\bf 269} (2007), 311--365.

\bibitem{Sch} E. Scheck, {\em Weyl laws for partially open quantum maps,}
Annales Henri Poincar\'e, {\bf 10} (2009), 714--747.


\bibitem{Sj} 
J. Sj\"ostrand,
{\em Eigenvalue distribution for non-self-adjoint operators on compact manifolds with small multiplicative random perturbations,} {\tt  arXiv:0809.4182}.



\bibitem{SjZw07}  J. Sj\"ostrand and M. Zworski,
{\em Elementary linear algebra for advanced spectral 
problems,} Ann. Inst. Fourier (Grenoble) {\bf 57} (2007), 2095-2141. 

\bibitem{SjZw04} J. Sj\"ostrand and M. Zworski,
{\em Fractal upper bounds on the density of semiclassical resonances,}
Duke Math. J. {\bf 137} (2007), 381--459.

\bibitem{Zw} M. Zworski, {\em Numerical linear algebra and solvability of partial differential equations}, Comm. Math. Phys. {\bf 229} (2002), 293--307.

\end{thebibliography}
\end{document}